\documentclass[]{article}
\usepackage{graphicx}
\usepackage{epstopdf}
\usepackage{amsfonts}
\usepackage{amsmath}
\usepackage{amssymb}
\usepackage{fancyhdr}
\usepackage{titlesec}
\usepackage{indentfirst}
\usepackage{booktabs}
\usepackage{verbatim}
\usepackage{color}
\usepackage{amsthm}
\usepackage{subfigure}

\usepackage[page,header]{appendix}
\usepackage{titletoc}

\numberwithin{equation}{section}
\newtheorem{theorem}{Theorem}[section]

\newtheorem{remark}[theorem]{Remark}

\newcommand{\rd}{\,\mathrm{d}}

\begin{document}
\title{An efficient dynamical low-rank algorithm \\for the Boltzmann--BGK equation \\close to the compressible viscous flow regime\footnote{
    JH's research was supported in part by NSF CAREER grant DMS-1654152 and NSF CDS\&E grant CBET-1854829.
    LY's research was supported in part by NSF DMS grant DMS-1818449.
  }
}

	\author{
	Lukas Einkemmer\footnote{Department of Mathematics, University of Innsbruck, Innsbruck, A-6020, Austria (lukas.einkemmer@uibk.ac.at).}, \
	Jingwei Hu\footnote{Department of Mathematics, Purdue University, West Lafayette, IN 47907, USA (jingweihu@purdue.edu).}, \  
	and \ Lexing Ying\footnote{Department of Mathematics, Stanford University, Stanford, CA 94305, USA (lexing@stanford.edu).}}    
	\maketitle
	
\begin{abstract}
    It has recently been demonstrated that dynamical low-rank algorithms can provide robust and efficient approximation to a range of kinetic equations. This is true especially if the solution is close to some asymptotic limit where it is known that the solution is low-rank. A particularly interesting case is the fluid dynamic limit that is commonly obtained in the limit of small Knudsen number. However, in this case the Maxwellian which describes the corresponding equilibrium distribution is not necessarily low-rank; because of this, the methods known in the literature are only applicable to the weakly compressible case. In this paper, we propose an efficient dynamical low-rank integrator that can capture the fluid limit -- the Navier-Stokes equations -- of the Boltzmann-BGK model even in the compressible regime. This is accomplished by writing the solution as $f=Mg$, where $M$ is the Maxwellian and the low-rank approximation is only applied to $g$. To efficiently implement this decomposition within a low-rank framework requires,  \textcolor{black}{in the isothermal case,} that certain coefficients are evaluated using convolutions, for which fast algorithms are known. Using the proposed decomposition also has the advantage that the rank required to obtain accurate results is significantly reduced compared to the previous state of the art. We demonstrate this by performing a number of numerical experiments and also show that our method is able to capture sharp gradients/shock waves.
\end{abstract}

{\small 

}

\section{Introduction}
\label{sec:intro}

The Navier-Stokes (NS) equations are widely used in many scientific and engineering disciplines and its fundamental importance cannot be underestimated. Nevertheless, it is increasingly realized that the NS equations may not provide an accurate description when the underlying system is not dense enough to be in thermodynamic equilibrium (i.e.,~when there is a significant deviation from the Maxwellian distribution). This is especially true when the mean free path is comparable to the characteristic length, a situation that occurs, e.g., in the atmosphere at high altitude or in micro-devices. To accurately describe such rarefied or transitional flows, the kinetic Boltzmann equation should be used \cite{Cercignani}. However, this comes at the price of solving a six-dimensional equation in phase space. Although significant progress has been made in the past decades to solve Boltzmann type kinetic equations more accurately and efficiently (see for instance \cite{Pareschi}), the high-dimensional nature of such equations still poses a challenging problem in scientific computing. It is therefore the goal of this paper to introduce an efficient method for the kinetic BGK equation in a wide range of flow regimes whose computational cost is comparable to solving the macroscopic fluid equations. 

The BGK equation, initially proposed by Bhatnagar, Gross, and Krook \cite{BGK54}, is one of the widely used kinetic models as it is much simpler than the full Boltzmann equation yet possesses most of its key properties. In dimensionless form, the equation reads 
\begin{equation} \label{BGK}
\partial_{t}f+v\cdot\nabla_{x}f=\frac{\nu}{\varepsilon}(M-f),\quad t>0,\quad x\in\Omega\subset\mathbb{R}^{d_x},\quad v\in\mathbb{R}^{d_v}, 
\end{equation}
where $f=f(t,x,v)$ is the one-particle probability density function depending on time $t$, position $x$, and velocity $v$. The so-called Maxwellian $M$ is defined as
\begin{equation} \label{Max}
M(t,x,v)=\frac{\rho(t,x)}{(2\pi T(t,x))^{d_v/2}}\exp\left(-\frac{|v-u(t,x)|^{2}}{2T(t,x)}\right)
\end{equation}
with the density $\rho$, bulk velocity $u$, and temperature $T$ given by the moments of $f$: 
\begin{equation} \label{moments}
\rho=\int_{\mathbb{R}^{d_v}}f\rd{v},\quad u=\frac{1}{\rho}\int_{\mathbb{R}^{d_v}}vf\rd{v},\quad T=\frac{1}{d_v\rho}\int_{\mathbb{R}^{d_v}}|v-u|^{2}f\rd{v}.
\end{equation}
The viscosity $\nu$ is a function of $\rho$ and $T$; typically, one assumes $\nu=\rho T^{1-\omega}$ with $0.5\leq \omega \leq 1$. Finally, $\varepsilon$ is the Knudsen number, defined as the ratio of the mean free path and the chosen characteristic length scale. The value of $\varepsilon$ indicates the flow regime \cite{Struchtrup} (see Section~\ref{sec:limit} for more details): 1) $\varepsilon \rightarrow 0$, Euler regime (the flow is well described by the compressible Euler equations; 2) $0<\varepsilon \lesssim 0.01$, NS regime (the flow is well described by the compressible NS equations); and 3) $0.01  \lesssim \varepsilon  \lesssim 1$, transition regime (the NS equations fail, and one has to resort to extended macroscopic models or the original kinetic equation (\ref{BGK})). 

In practical applications we often have $d_x=d_v=3$, hence the computational cost of solving (\ref{BGK}) would be at least $\mathcal{O}(N^6)$, where $N$ is the number of mesh points in each dimension of the physical and velocity space. To reduce this cost, we are going to adopt the basic framework of so-called dynamical low-rank projection, which has a long history in quantum dynamics; see, for example, \cite{meyer90tmc,meyer09mqd} for the MCTDH approach to molecular quantum dynamics in the chemical physics literature and \cite{Lubich2008,lubich15tii,Conte2010} for a mathematical point of view. In a general mathematical setting, the dynamical low-rank approximation has been studied in \cite{Koch2007,Koch2010,lubich13dab,arnold2014approximation}. A major algorithmic advance for the time integration was achieved with the projector-splitting method first proposed in \cite{Lubich2014}. In contrast to standard time-stepping methods, the projector-splitting method is robust to the presence of small singular values in the low-rank approximation~\cite{Kieri2016}. Recently, this method has been applied to several kinetic equations including the Vlasov equation \cite{El18, EL18_cons, piazzola}, the BGK equation (in the weakly compressible regime) \cite{E18}, and radiative transport equations \cite{peng2019low,DEL19, EHW20}.

Though the details vary for each of the aforementioned equations, the common idea is to seek a low-rank approximation of the unknown function $f$ as
\begin{equation} \label{lowrank1}
f(t,x,v)\approx \sum_{i,j=1}^r X_i(t,x)S_{ij}(t)V_j(t,v),
\end{equation}
where $r$ is the rank, $\{X_i\}$ and $\{V_j\}$ are the orthonormal bases in the physical space and velocity space, respectively. The dynamics of $X_i$, $S_{ij}$ and $V_j$ are determined by projecting the equation onto the tangent space of the low-rank solution manifold. Upon a further operator splitting, the original six-dimensional problem can be reduced to a few three-dimensional problems. It is also clear that the efficiency of this low-rank method depends on the intrinsic rank of the underlying problem. If the solution $f$ has no low-rank structure, a large value of $r$ would be needed in (\ref{lowrank1}) and hence the gain of the method would be marginal. For a number of linear and weakly nonlinear problems it is known that the solution is low-rank in some asymptotic limit. This is the case, for instance, in the radiative transfer/linear transport equations wherein one obtains a rank-$1$ solution in the diffusion limit. The asymptotic-preserving (AP) dynamical low-rank methods that can capture this behavior have been considered in \cite{DEL19,EHW20}. A similar approach was undertaken in \cite{E18} for the BGK equation, where the method is limited to the weakly compressible flows since the Maxwellian \eqref{Max} is not necessarily low-rank.

In this work we introduce a robust dynamical low-rank integrator for the BGK equation (\ref{BGK}) that is able to capture the fluid limit -- the NS equations -- in the compressible (isothermal) case. We do this by writing the solution as $f = M g$, where only $g$ is subjected to a low-rank approximation. Note that this is a drastically different decomposition compared to (\ref{lowrank1}) due to the  nonlinearity of $M$. To efficiently implement this decomposition within the low-rank framework requires the evaluation of certain integrals involving the Maxwellian. If this is done naively the numerical method would suffer from excessive computational cost that scales as $\mathcal{O}(N^6)$. \textcolor{black}{We demonstrate that, in the isothermal case, this difficulty can be} overcome by interpreting the corresponding integrals as convolutions and using a FFT based fast convolution algorithm.

The proposed approach has several advantages. First, from the Chapman-Enskog expansion \cite{BGL91} we can guarantee that $g$ is low-rank for small values of $\varepsilon$ (in both the Euler and NS regimes) for compressible problems. Second, for weakly compressible problems, the rank $r$ needed to obtain an accurate approximation is much smaller than that in \cite{E18} (because we only assume that $g=f/M$ is low-rank and do not need to perform the expansion of $M$ in small velocities as was done in \cite{E18}). Third, the fluid moments, such as density $\rho$ and momentum $\rho u$, are integrated as part of the numerical scheme. Thus, using an appropriate discretization, as will be done here, conservation of mass and momentum is ensured. The algorithm is also able to handle problems with sharp gradients/shock waves. In addition, since the proposed low-rank algorithm automatically takes care of regions close to the fluid regime (the scheme becomes a consistent discretization to the NS equations when $\varepsilon$ is small), it can be used to perform simulations in a wide range of flow regimes without the need to explicitly couple a kinetic solver with a fluid one (which is potentially complicated in a hybrid approach; see, e.g., \cite{filbet2015,crouseilles2005hybrid}). In that context the present scheme would be used in a similar way as the Hermite expansion in \cite{Vencels2015}, with the added advantage that the computational cost of the kinetic solver is reduced dramatically.

The rest of this paper is organized as follows. In Section~\ref{sec:fluid-limit} we describe the fluid dynamic limits of the BGK equation. We then motivate the low-rank approximation used (Section~\ref{sec:low-rank}) and derive the corresponding numerical algorithm in a semi-discrete setting (Section~\ref{sec:low-rank2}). The full discretization is discussed in Section~\ref{sec:fully-discrete}. A formal analysis is performed in Section~\ref{sec:ap} to show that our numerical scheme can indeed capture the NS limit when $\varepsilon$ is small. Finally we present the results of a range of numerical experiments in Section~\ref{sec:numerical}.

\section{Fluid limits of the BGK equation \label{sec:fluid-limit}}
\label{sec:limit}

In this section, we derive the fluid dynamic limits of the BGK equation (\ref{BGK}) using the Chapman-Enskog expansion \cite{BGL91}. We will see that the compressible Euler equations can be obtained when $\varepsilon\rightarrow0$, while the compressible NS equations can be obtained when retaining $\mathcal{O}(\varepsilon)$ information. Although the procedure is quite standard in the kinetic literature, we have decided to include it here because our following low-rank approximation is strongly motivated by this analysis.

When $\varepsilon$ is small, from (\ref{BGK}), formally one has
\begin{equation}
f=M+\mathcal{O}(\varepsilon).
\end{equation}
Hence it makes sense to assume that
\begin{equation} \label{ansatz}
f=M+\varepsilon f_1,
\end{equation}
where $f_1$ is an $\mathcal{O}
(1)$ function. We then substitute (\ref{ansatz}) into (\ref{BGK}) to obtain
\begin{equation} \label{ff}
f_1=-\frac{1}{\nu}(\partial_tf+v\cdot \nabla_x f)=-\frac{1}{\nu}(\partial_tM+v\cdot \nabla_x M)+\mathcal{O}(\varepsilon).
\end{equation}

On the other hand, taking the first $d_v+2$ moments of (\ref{BGK}), i.e.~multiplying (\ref{BGK}) by $\phi(v):=(1,v,|v|^{2}/2)^{T}$ and integrating in $v$, yields
\begin{equation} \label{UU}
\partial_{t}\langle f\phi \rangle_v +\nabla_{x}\cdot\langle v\phi f\rangle_{v}=0, \quad \langle \ \cdot \ \rangle_v:= \int_{\mathbb{R}^{d_v}} \cdot \ \rd{v}.
\end{equation}
By substituting (\ref{ansatz}) into (\ref{UU}) we obtain
\begin{equation} \label{UU1}
\partial_{t}\langle f\phi \rangle_v+\nabla_{x}\cdot\langle v\phi M \rangle_{v}=-\varepsilon \nabla_{x}\cdot\langle v\phi f_1\rangle_{v}.
\end{equation}
Using the macroscopic quantities defined in (\ref{moments}) and
\begin{equation}
\mathbb{P}_1:=-\int_{\mathbb{R}^{d_v}} (v-u)\otimes (v-u)f_1\rd{v}, \quad q_1:=-\frac{1}{2}\int_{\mathbb{R}^{d_v}} (v-u)|v-u|^2f_1\rd{v},
\end{equation}
we can write (\ref{UU1}) as
\begin{align}  \label{NS1}
\left\{
\begin{array}{l}
\partial_t \rho+\nabla_x\cdot (\rho u)=0,\\[8pt]
\partial_t (\rho u)+\nabla_x\cdot (\rho u\otimes u +p\text{Id})=\varepsilon\nabla_x \cdot \mathbb{P}_1,\\[8pt]
\partial_t E+\nabla_x\cdot ((E+p)u)=\varepsilon \nabla_x \cdot (\mathbb{P}_1u+q_1),
\end{array}\right.
\end{align}
where $p=\rho T$ is the pressure, Id is the identity matrix, and $E=\frac{d_v}{2}\rho T+\frac{1}{2}\rho u^2$ is the total energy. Neglecting the $\mathcal{O}(\varepsilon)$ terms, (\ref{NS1}) become the compressible Euler equations.

To further derive the right hand side of (\ref{NS1}), we first note that for the Maxwellian function defined in (\ref{Max}) one has
\begin{equation} \label{MM2}
\frac{1}{M}(\partial_t M+v \cdot \nabla_x M)=\frac{1}{\rho} (\partial_t \rho+v \cdot \nabla_x \rho)+\frac{(v-u)}{T}\cdot (\partial_t u+v \cdot \nabla_x u)+\left ( \frac{|v-u|^2}{2T^2}-\frac{d_v}{2T}\right) (\partial_t T+v \cdot \nabla_x T).
\end{equation}
Since we already knew that $\rho$, $u$ and $T$ satisfy the compressible Euler equations to the leading order, we can use the spatial derivatives of these functions to replace the time derivatives. This yields 
\begin{equation} \label{MM1}
\frac{1}{M}(\partial_t M+v \cdot \nabla_xM)= \left(\frac{|v-u|^2}{2T}-\frac{d_v+2}{2}\right)\frac{(v-u)\cdot \nabla_xT}{T}+\left( \frac{(v-u)\otimes (v-u)}{T}-\frac{|v-u|^2}{d_vT} \text{Id}\right):\nabla_xu+\mathcal{O}(\varepsilon),
\end{equation}
where
\begin{equation} \label{sigmau}
\sigma(u):=\nabla_xu+(\nabla_x u)^T- \frac{2}{d_v}(\nabla_x\cdot u) \ \text{Id}.
\end{equation}
Note that $\nabla_x u$ is a matrix with $ij$-th component given by $\partial_{x_j} u_i$ and the operation $:$ between two matrices is defined as $A:B=\sum_{ij}a_{ij}b_{ij}$.

Plugging (\ref{MM1}) into (\ref{ff}) gives
\begin{equation} \label{ff1}
f_1=-\frac{M}{\nu}\left[\left(\frac{|v-u|^2}{2T}-\frac{d_v+2}{2}\right)\frac{(v-u)\cdot \nabla_xT}{T}+\left( \frac{(v-u)\otimes (v-u)}{T}-\frac{|v-u|^2}{d_vT} \text{Id}\right):\nabla_xu\right]+\mathcal{O}(\varepsilon).
\end{equation}
One can verify by direct integration that
\begin{equation} 
\begin{split}
\mathbb{P}_1&=\frac{1}{\nu}\rho T\sigma(u)+\mathcal{O}(\varepsilon)=T^\omega \sigma(u)+\mathcal{O}(\varepsilon),\\
q_1&=\frac{1}{\nu}\frac{d_v+2}{2}\rho T\nabla_xT+\mathcal{O}(\varepsilon)=\frac{d_v+2}{2}T^\omega \nabla_xT+\mathcal{O}(\varepsilon).
\end{split}
\end{equation}
Substituting the above $\mathbb{P}_1$ and $q_1$ into (\ref{NS1}), we finally obtain
\begin{align} \label{NS}
\left\{
\begin{array}{l}
\partial_t \rho+\nabla_x\cdot (\rho u)=0,\\[8pt]
\partial_t (\rho u)+\nabla_x\cdot (\rho u\otimes u +p\text{Id})=\varepsilon\nabla_x \cdot (\mu \sigma(u)) +\mathcal{O}(\varepsilon^2),\\[8pt]
\partial_t E+\nabla_x\cdot ((E+p)u)=\varepsilon \nabla_x \cdot (\mu \sigma(u)u+\kappa \nabla_x T)+\mathcal{O}(\varepsilon^2),
\end{array}\right.
\end{align}
where $\mu:=T^{\omega}$ and $\kappa:=\frac{d_v+2}{2}\mu$ are, respectively, the coefficients of viscosity and heat conductivity. Ignoring $\mathcal{O}(\varepsilon^2)$ terms, (\ref{NS}) are the compressible Navier-Stokes equations. 

\begin{remark}
It is well-known that the BGK model (\ref{BGK}) cannot produce the correct Prandtl number ($\text{Pr}:=\frac{d_v+2}{2}\frac{\mu}{\kappa}\equiv 1$ in the BGK model while for most gases $\text{Pr}<1$). To correct this defect, various models have been proposed such as the ES-BGK model \cite{ALPP00} or the Shakhov model \cite{Shakhov68}. Since our main focus in this paper are isothermal flows, this defect is irrelevant. 
\end{remark}

\section{The low-rank approximation \label{sec:low-rank}}

From the previous section we can see that when $\varepsilon$ is small, $f$ can be expanded as (see  (\ref{ansatz}) and (\ref{ff1})):
\begin{equation}
f=M-\varepsilon \frac{M}{\nu}\left[\left(\frac{|v-u|^2}{2T}-\frac{d_v+2}{2}\right)\frac{(v-u)\cdot \nabla_xT}{T}+\left( \frac{(v-u)\otimes (v-u)}{T}-\frac{|v-u|^2}{d_vT} \text{Id}\right):\nabla_xu\right]+\mathcal{O}(\varepsilon^2).
\end{equation}
Note that the Maxwellian $M$ is a normal distribution with (time and spatially dependent) mean $u(t,x)$ and variance $T(t,x)$, so $M$ is generally not a low-rank or separable function in $x$ and $v$. Because of this, we do not directly use the low-rank approximation (\ref{lowrank1}) for $f$. Instead, we consider
\begin{equation}
f=Mg,
\end{equation}
then
\begin{equation}
g =1-\varepsilon\frac{1}{\nu}\left[\left(\frac{|v-u|^{2}}{2T}-\frac{d_v+2}{2}\right)\frac{(v-u)\cdot\nabla_{x}T}{T}+\left(\frac{(v-u)\otimes(v-u)}{T}-\frac{|v-u|^{2}}{d_vT}\text{Id}\right):\nabla_{x}u\right]+\mathcal{O}(\varepsilon^{2}),\label{gg}
\end{equation}
which is a low-rank function in $x$ and $v$ even at $\mathcal{O}(\varepsilon)$. \textcolor{black}{This can be most easily seen by recognizing that the $\mathcal{O}(\varepsilon)$ term can be written as a sum of products of functions that depend only on $x$ (such as $u$, $1/T$, and $\nabla_x u$) and functions that only depend on $v$ (such as $v$, $v^2$ and $v \otimes v$).}

This motivates us to seek an expansion of $g$ in the following form: 
\begin{equation}
g=\sum_{i,j=1}^{r}X_{i}(t,x)S_{ij}(t)V_{j}(t,v),\label{lowrank}
\end{equation}
where $\{X_{i}\}$ and $\{V_{j}\}$ are the orthonormal bases in $x$ and $v$, respectively.
\textcolor{black}{Before proceeding, let us note that a multiplicative deviation from equilibrium to improve compression has been used in the context of sparse grids and the Vlasov equation \cite{kormann2016sparse}, albeit in the somewhat simpler setting where $u=0$.}

To track the evolution of $M$, one can consider the moment equation (\ref{UU}) which we rewrite using a slightly different notation:
\begin{equation}
\partial_{t}U+\nabla_{x}\cdot\langle v\phi Mg\rangle_{v}=0,\label{MM}
\end{equation}
where $U:=(\rho,\rho u,E)^{T}$ and $\phi(v):=(1,v,|v|^{2}/2)^{T}$. 

To track the dynamics of $X_{i}$, $S_{ij}$ and $V_{j}$, we proceed as follows. Substituting $f=Mg$ into equation (\ref{BGK}) yields
\begin{equation}
\partial_{t}g=-v\cdot\nabla_{x}g-\frac{1}{M}(\partial_{t}M+v\cdot\nabla_{x}M)g+\frac{\nu}{\varepsilon}(1-g):=h.\label{hh}
\end{equation}
We now apply the projector-splitting based dynamical low-rank algorithm \cite{Lubich2014} to \eqref{hh}. That is, constraining $g$ to lie on the low-rank manifold:
\begin{equation}
\partial_{t}g=\sum_{j}\langle V_{j},h\rangle_{v}V_{j}-\sum_{i,j}X_{i}\langle X_{i}V_{j},h\rangle_{x,v}V_{j}+\sum_{i}X_{i}\langle X_{i},h\rangle_{x}
\end{equation}
and then performing the operator splitting yields the following numerical algorithm:

-- Let $K_{j}=\sum_{i}X_{i}S_{ij}$, then $g=\sum_{j}K_{j}V_{j}$. Update $K_{j}$ by solving
\begin{equation}
\partial_{t}K_{j}=\langle V_{j},h\rangle_{v},\label{eq:evol-K}
\end{equation}
then perform an orthonormalization of $K_{j}$ to generate new $X_{i}$ and $S_{ij}$ (using a QR decomposition).

-- Update $S_{ij}$ by solving
\begin{equation}
\partial_{t}S_{ij}=-\langle X_{i}V_{j},h\rangle_{x,v}.\label{eq:evol-S}
\end{equation}

-- Let $L_{i}=\sum_{j}S_{ij}V_{j}$, then $g=\sum_{i}X_{i}L_{i}$. Update $L_{i}$ by solving
\begin{equation}
\partial_{t}L_{i}=\langle X_{i},h\rangle_{x},\label{eq:evol-L}
\end{equation}
then perform an orthonormalization on $L_{i}$ to generate new $S_{ij}$ and $V_{j}$ (using a QR decomposition).

For the derivation of the above algorithm we refer the readers to \cite{El18} and the references contained therein. The algorithm described is robust in the sense that over approximation, i.e.~small singular values in $S$, do not cause any numerical issues. Furthermore, the order of the subflows \eqref{eq:evol-K}-\eqref{eq:evol-L} can be chosen freely. Any order results in a numerical method that is first order accurate (in time). For higher order methods this can be exploited to lower the overall computational cost (see, e.g., \cite{El18}).

\section{The dynamical low-rank algorithm\label{sec:low-rank2}}

In this section we present the proposed algorithm in detail. This includes the low-rank integrator for $g$ as well as the necessary numerical scheme to advance the moments in time. The present section describes a semi-discrete version where time is discretized but space is left continuous. We will consider possible spatial discretizations in more details in Section~\ref{sec:fully-discrete}.

From this point onward, we will assume that the temperature is constant in order to simplify the presentation. Without loss of generality we can then choose $T \equiv 1$. Assuming that the ideal gas law holds true, this simply means that we measure the velocity $u$ in fractions of the speed of sound. \textcolor{black}{The rationale for this assumption lies in the fact that isothermal flows represent an interesting nontrivial compressible case that we propose to efficiently treat using fast convolution algorithms.}  We will discuss the non-constant temperature case in more detail in Remark~\ref{remark:temp} given at the end of this section.

First, if $T=1$ in (\ref{BGK}), we have $\nu=\rho$ and the Maxwellian becomes 
\begin{equation}
M=\frac{\rho}{(2\pi)^{d_v/2}}\exp\left(-\frac{|v-u|^{2}}{2}\right).
\end{equation}
The equations of motions for the moments (\ref{MM}) can then be written as
\begin{align}
  & \partial_{t}\rho=-\nabla_{x}\cdot\left( \sum_{i,j} X_{i}S_{ij}\langle vV_{j}M\rangle_{v}\right):=I_{1},\label{rho}\\
  & \partial_{t}(\rho u)=-\nabla_{x}\cdot\left(\sum_{i,j} X_{i}S_{ij}\langle v\otimes vV_{j}M\rangle_{v}\right):=I_{2}, \label{rhou}
\end{align}
from which we can easily derive the following expression
\begin{equation}\label{u}
\partial_{t}u=\frac{1}{\rho}(I_{2}-I_{1}u).
\end{equation}
Furthermore, $h$ in (\ref{hh}) becomes 
\begin{equation}
h=-v\cdot\nabla_{x}g-\mathcal{M}g+\frac{\rho}{\varepsilon}(1-g), \label{hh1}
\end{equation}
where
\begin{equation} \label{MM123}
\begin{split}
\mathcal{M}:&=\frac{1}{M}(\partial_{t}M+v\cdot\nabla_{x}M)=\frac{1}{\rho}(\partial_{t}\rho+v\cdot\nabla_{x}\rho)+(v-u)\cdot(\partial_{t}u+v\cdot\nabla_{x}u)\\
&=\frac{1}{\rho}\left(I_{1}-u\cdot(I_{2}-I_{1}u)\right)+v\cdot\left(\frac{1}{\rho}\left(\nabla_{x}\rho+I_{2}-I_{1}u\right)-\frac{1}{2}\nabla_{x}(u^{2})\right)+(v\otimes v):\nabla_{x}u\\
:&=\mathcal{M}_{1}+v\cdot\mathcal{M}_{2}+(v\otimes v):\mathcal{M}_{3},
\end{split}
\end{equation}
where
\begin{equation} \label{M123}
\mathcal{M}_{1}  =\frac{1}{\rho}(I_{1}-u\cdot(I_{2}-I_{1}u)), \quad \mathcal{M}_{2}  =\frac{1}{\rho}(\nabla_{x}\rho+I_{2}-I_{1}u)-\frac{1}{2}\nabla_{x}(u^2), \quad \mathcal{M}_{3}  =\nabla_{x}u.
\end{equation}

We now plug (\ref{hh1}) into (\ref{eq:evol-K})  to obtain 
\begin{equation} \label{KK}
\begin{split}
\partial_{t}K_{j} & =\langle V_{j},h\rangle_{v}\\
 & =-\sum_{kl}\langle V_{j},v\cdot(\nabla_{x}X_{k})S_{kl}V_{l}\rangle_{v}-\sum_{kl}X_{k}S_{kl}\langle V_{j},V_{l}\mathcal{M}\rangle_{v}+\frac{\rho}{\varepsilon}\left(\langle V_{j}\rangle_{v}-\sum_{kl}X_{k}S_{kl}\langle V_{j},V_{l}\rangle_{v}\right)\\
 & =-\sum_{l}(\nabla_{x}K_{l})\cdot\langle vV_{j}V_{l}\rangle_{v}-\sum_{l}K_{l}\langle V_{j}V_{l}\mathcal{M}\rangle_{v}+\frac{\rho}{\varepsilon}\left(\langle V_{j}\rangle_{v}-K_{j}\right),
 \end{split}
\end{equation}
where
\begin{equation}
\langle V_{j}V_{l}\mathcal{M}\rangle_{v}=\delta_{jl}\mathcal{M}_{1}+\langle vV_{j}V_{l}\rangle_{v}\cdot\mathcal{M}_{2}+\langle v\otimes v V_{j}V_{l}\rangle_{v}:\mathcal{M}_{3}.
\end{equation}
Plugging (\ref{hh1}) into (\ref{eq:evol-S}) gives
\begin{equation} \label{SS}
\begin{split}
\partial_{t}S_{ij} & =-\langle X_{i}V_{j},h\rangle_{xv}\\
 & =\sum_{kl}(\langle X_{i}\nabla_{x}X_{k}\rangle_{x}\cdot\langle vV_{j}V_{l}\rangle_{v})S_{kl}+\sum_{kl}S_{kl}\langle X_{i}X_{k}V_{j}V_{l}\mathcal{M}\rangle_{xv}-\frac{1}{\varepsilon}\left(\langle \rho X_{i}\rangle_{x}\langle V_{j}\rangle_{v}-\sum_{k}\langle \rho X_k X_i\rangle_x S_{kj}\right),
\end{split}
\end{equation}
where
\begin{equation} 
\begin{split}
\langle X_{i}X_{k}V_{j}V_{l}\mathcal{M}\rangle_{xv} & =\langle X_{i},X_{k}\langle V_{j}V_{l}\mathcal{M}\rangle_{v}\rangle_{x}\\
 & =\delta_{jl}\langle X_{i}X_{k}\mathcal{M}_{1}\rangle_{x}+\langle vV_{j}V_{l}\rangle_v\cdot\langle X_{i}X_{k}\mathcal{M}_{2}\rangle_{x}+\langle v\otimes vV_{j}V_{l}\rangle_v\colon\langle X_{i}X_{k}\mathcal{M}_{3}\rangle_{x}.
 \end{split}
\end{equation}
Finally plugging (\ref{hh1}) into (\ref{eq:evol-L}) gives
\begin{equation}  \label{LL}
\begin{split}
\partial_{t}L_{i} & =\langle X_{i},h\rangle_{x}\\
 & =-\sum_{kl}\langle X_{i},\nabla_{x}X_{k}\rangle_{x}\cdot vS_{kl}V_{l}-\sum_{kl}\langle X_{i},X_{k}\mathcal{M}\rangle_{x}S_{kl}V_{l}+\frac{1}{\varepsilon}\left(\langle \rho X_{i}\rangle_{x}-\sum_k \langle \rho X_k X_i\rangle L_k \right)\\
 & =-\sum_{k}(\langle X_{i}\nabla_{x}X_{k}\rangle_{x}\cdot v)L_{k}-\sum_{k}\langle X_{i}X_{k}\mathcal{M}\rangle_{x}L_{k}+\frac{1}{\varepsilon}\left(\langle \rho X_{i}\rangle_{x}- \sum_k \langle \rho X_k X_i\rangle L_k\right),
 \end{split}
\end{equation}
where
\begin{equation}
\langle X_{i}X_{k}\mathcal{M}\rangle_{x}=\langle X_{i}X_{k}\mathcal{M}_{1}\rangle_{x}+v\cdot\langle X_{i}X_{k}\mathcal{M}_{2}\rangle_{x}+(v\otimes v):\langle X_{i}X_{k}\mathcal{M}_{3}\rangle_x.
\end{equation}

\subsection{Time discretization \label{sec:first-order}}

We are now ready to present the semi-discrete \textcolor{black}{(discrete in time and continuous in space)} dynamical low-rank integrator. The main ingredients of this scheme are the low-rank approximation described in the previous section (i.e.~$f=Mg$ and the low-rank expansion is applied to $g$), using fast convolution algorithms to evaluate the integrals in \eqref{eq:upmom1}-\eqref{eq:upmom2} below, and using an IMEX scheme to capture the fluid limit. We will also provide the computational cost for each substep to highlight the efficiency of the proposed algorithm compared to full grid simulation that scale as $\mathcal{O}(N^{d_x+d_v})$.

Suppose that at time step $t^n$, $(\rho^n, u^n, X_i^n,V_j^n,S_{ij}^n)$ are available. In order to
obtain the solution $(\rho^{n+1}, u^{n+1}, X_i^{n+1},V_j^{n+1},S_{ij}^{n+1})$ at $t^{n+1}$, we
proceed as follows:

\vspace{0.1in}
\noindent \textbf{Step 1: Update $\rho^{n+1}$ and $u^{n+1}$.}
\begin{enumerate}
\item Compute the terms
\begin{align}
  & \langle vV_{j}^{n}M^{n}\rangle_{v} =\frac{\rho^{n}(x)}{(2\pi)^{d_v/2}}\left\langle vV_{j}^{n}(v)\exp\left(-\frac{|v-u^n(x)|^{2}}{2}\right)\right\rangle_{v}, \label{eq:upmom1}\\
  &\langle v\otimes vV_{j}^{n}M^{n}\rangle_{v}  =\frac{\rho^{n}(x)}{(2\pi)^{d_v/2}}\left\langle v\otimes vV_{j}^{n}(v)\exp\left(-\frac{|v-u^n(x)|^{2}}{2}\right)\right\rangle_{v}.\label{eq:upmom2}
\end{align}
Note that these two quantities can be written as convolutions and hence are computed as follows:

-- Compute the convolutions
\begin{align*}
  &g_{j}^{1} =(v\mapsto vV_{j}^{n})\ast(v\mapsto\exp(-v^{2}/2)),\\
  &g_{j}^{2} =(v\mapsto(v\otimes v)V_{j}^n)\ast(v\mapsto\exp(-v^{2}/2)),
\end{align*}
using an FFT. \textbf{Cost: $\mathcal{O}(rN^{d_v}\log N^{d_v})$. }\\
-- Evaluate
\begin{align*}
&\langle vV_{j}^{n}M^{n}\rangle_{v}(x) =\frac{\rho^{n}(x)}{(2\pi)^{d_v/2}}g_{j}^{1}(u^n(x)),\\
&\langle(v\otimes v)V_{j}^{n}M^{n}\rangle_{v}(x)  =\frac{\rho^{n}(x)}{(2\pi)^{d_v/2}}g_{j}^{2}(u^n(x)),
\end{align*}
using cubic splines (or any other interpolation). \textcolor{black}{This is done by taking the equidistant output obtained from the FFT and constructing the corresponding interpolant.} \textbf{Cost: $\mathcal{O}(rN^{d_{x}})$.}

Let us emphasize that using the described procedure is essential as computing \eqref{eq:upmom1}-\eqref{eq:upmom2} in a naive way would incur a computational cost of $\mathcal{O}(r N^{d_x+d_v})$, which is clearly prohibitive.

\item Compute 
\begin{align} 
 &   I_{1}^{n}  =-\nabla_{x}\cdot\left(\sum_{ij}X_{i}^{n}S_{ij}^{n}\langle vV_{j}^{n}M^{n}\rangle_{v}\right), \label{eq:I1n} \\
& I_{2}^{n}  =-\nabla_{x}\cdot\left(\sum_{ij}X_{i}^{n}S_{ij}^{n}\langle v\otimes vV_{j}^{n}M^{n}\rangle_{v}\right).\label{eq:I2n} 
\end{align}
\textbf{Cost: $\mathcal{O}(r^{2}N^{d_{x}})$.}

\item Perform a forward Euler step to solve (\ref{rho}) and (\ref{u}):
\begin{align}
 &   \rho^{n+1}  =\rho^{n}+\Delta tI_{1}^{n}, \label{eq:rhon+1} \\
 &   u^{n+1}  =u^{n}+\frac{1}{\rho^{n}}(I_{2}^{n}-I_{1}^{n}u^{n}). \label{eq:un+1}
\end{align}
\textbf{Cost: $\mathcal{O}(N^{d_{x}}).$}
\end{enumerate}

\vspace{0.1in}
\noindent \textbf{Step 2: Update $X_i^{n+1}$, $V_j^{n+1}$, and $S_{ij}^{n+1}$.}

\vspace{0.1in}
\noindent \textbf{K step}
\begin{enumerate}
\item Compute
\begin{align*}
c_{jl}^{1}  =\langle vV_{j}^{n}V_{l}^{n}\rangle_{v}, \quad c_{jl}^{\star} =\langle v\otimes v V_{j}^{n}V_{l}^{n}\rangle_{v}, \quad \overline{V}_{j}  =\langle V_{j}^{n}\rangle_{v}.
\end{align*}
\textbf{Cost: $\mathcal{O}(r^{2}N^{d_v}).$}
\item Compute $\mathcal{M}_1$, $\mathcal{M}_2$, and $\mathcal{M}_3$ defined in (\ref{M123}) using $\rho^n$, $u^n$, $I_1^n$, and $I_2^n$. \textbf{Cost: $\mathcal{O}(N^{d_{x}}).$ }
\item Compute 
\[
c_{jl}^{2}=\delta_{jl}\mathcal{M}_{1}+c_{jl}^{1}\cdot\mathcal{M}_{2}+c_{jl}^{\star}:\mathcal{M}_3.
\]
\textbf{Cost: $\mathcal{O}(r^{2}N^{d_{x}}).$}
\item Perform a first order IMEX step to solve (\ref{KK}):
\begin{align}
&K_{j}^{n}  =\sum_{i}X_{i}^{n}S_{ij}^{n}, \nonumber \\
 &   K_{j}^{n+1} =\frac{1}{1+\Delta t \rho^n/\varepsilon}K_{j}^{n}-\frac{\Delta t}{1+\Delta t \rho^n/\varepsilon}\left[\sum_{l}c_{jl}^{1}\cdot(\nabla_{x}K_{l}^{n})+\sum_{l}c_{jl}^{2}K_{l}^{n}\right]+\frac{\Delta t \rho^n}{\varepsilon+\Delta t \rho^n}\overline{V}_{j}. \label{eq:Kn+1}
\end{align}
\textbf{Cost: $\mathcal{O}(r^{2}N^{d_{x}}).$} Note: The IMEX scheme is necessary to treat the stiff term and capture the asymptotic limit. Nevertheless, as can be seen, the update is completely explicit. For more details we refer the readers to \cite{EHW20}.
\item Compute a QR decomposition of $K_{j}^{n+1}$ to obtain $X_{i}^{n+1}$
and $S_{ij}^{1}$. \textbf{Cost: $\mathcal{O}(r^{2}N^{d_{x}}).$}
\end{enumerate}

\noindent \textbf{S step}
\begin{enumerate}
\item Compute
\begin{align}
    &d_{ik}^{1} =\langle X_{i}^{n+1}\nabla_{x}X_{k}^{n+1}\rangle_{x}, \quad d_{ik}^{\star} =\langle X_{i}^{n+1}X_{k}^{n+1}\mathcal{M}_{1}\rangle_{x}, \quad d_{ik}^{\star\star}  =\langle X_{i}^{n+1}X_{k}^{n+1}\mathcal{M}_{2}\rangle_{x}, \label{eq:d1} \\
&d_{ik}^{\star\star\star}  =\langle X_{i}^{n+1}X_{k}^{n+1}\mathcal{M}_{3}\rangle_{x},\quad \overline{X}_{i}  =\langle  \rho^nX_{i}^{n+1}\rangle_{x}, \quad R_{ik}  = \langle  \rho^n X_i^{n+1} X_k^{n+1}\rangle. \nonumber
\end{align}
\textbf{Cost: $\mathcal{O}(r^{2}N^{d_{x}}).$}
\item Compute
\[
d_{ik;jl}^{2}=\delta_{jl}d_{ik}^{\star}+c_{jl}^{1}\cdot d_{ik}^{\star\star}+c_{jl}^{\star}\colon d_{ik}^{\star\star\star}.
\]
\textbf{Cost: }$\mathcal{O}(r^{4}).$ Note: We could also compute $d_{ik;jl}^{2}=\langle X_{i}X_{k}c_{jl}^{2}\rangle_{x}$ but the cost would be $\mathcal{O}(r^{4}N^{d_{x}})$.
\item Perform a first order IMEX step to solve (\ref{SS}):
\[
    \sum_k \left(I - \frac{\Delta t}{\varepsilon}R\right)_{ik} S_{kj}^{2}=S_{ij}^{1} +\Delta t \left[\sum_{kl}(d_{ik}^{1}\cdot c_{jl}^{1})S_{kl}^{1}+\sum_{kl}d_{ik;jl}^{2}S_{kl}^{1}\right]-\frac{\Delta t}{\varepsilon}\overline{X}_{i}\overline{V}_{j}.
\]
\textbf{Cost: $\mathcal{O}(r^{4}).$}
\end{enumerate}

\noindent  \textbf{L step}
\begin{enumerate}
\item Perform a first order IMEX step to solve (\ref{LL}):
\begin{align*}
&L_{i}^{n}  =\sum_{j}S_{ij}^{2}V_{j}^{n},\\
&\sum_{j} \left(I+\frac{\Delta t}{\varepsilon}R\right)_{ij} L_{j}^{n+1}  = L_{i}^{n}-\Delta t \left[\sum_{k}(d_{ik}^{1}\cdot v)L_{k}^{n}+\sum_{k}\left(d_{ik}^{\star}+v\cdot d_{ik}^{\star\star}+(v\otimes v):d_{ik}^{\star\star\star}\right)L_{k}^{n}\right]+\frac{\Delta t}{\varepsilon}\overline{X}_{i}.
\end{align*}
\textbf{Cost: $\mathcal{O}(r^{2}N^{d_v}).$}
\item Compute a QR decomposition of $L_{i}^{n+1}$ to obtain $V_{i}^{n+1}$ and $S_{ij}^{n+1}$. \textbf{Cost: $\mathcal{O}(r^{2}N^{d_v}).$}
\end{enumerate}

\begin{remark}
The scheme described is first order in time. It can be generalized to second (or higher) order as described in \cite{EHW20}. For simplicity we only consider the first order method in this paper.
\end{remark}

\begin{remark} \label{rem:reynolds}
It is worth mentioning that when $T=1$, the limiting system (\ref{NS}) become the isothermal NS equations:
\begin{align}
\left\{
\begin{array}{l}
\displaystyle \partial_t \rho+\nabla_x\cdot (\rho u)=0,\\[8pt]
\displaystyle \partial_t (\rho u)+\nabla_x\cdot (\rho u\otimes u +\rho\text{Id})=\varepsilon \nabla_x \cdot \sigma(u),
\end{array}\right.
\end{align}
where $\sigma(u)$ is the same as before and $1/\varepsilon$ plays the same role as the Reynolds number.
\end{remark}
 
\begin{remark} \label{remark:temp}
We have made the simplifying assumption $T=1$. \textcolor{black}{The algorithm can be extended easily to the case of
piecewise constant temperature $T^n(x)$ with two changes. First, besides the equations
of motions \eqref{rho} and \eqref{rhou}, there will be an additional equation of motion for
$E(t,x)$. Second, the computation of \eqref{eq:upmom1} or \eqref{eq:upmom2} can still be accelerated by
computing an FFT for each of the constant values of $T^n(x)$, followed by interpolating at $u^n(x)$
for each $x$ with that constant temperature value. For more general temperature distribution
$T^n(x)$, the situation is less clear and we leave this as future work.}

\end{remark}

\begin{remark}
    \textcolor{black}{The proposed dynamical low-rank integrator is currently the only such scheme that can be efficiently used in the compressible case. Nevertheless, it is instructive to compare its performance in the weakly compressible case with the method developed in \cite{E18}. The computational and memory complexity for both methods is, up to a logarithm, the same; namely, $\mathcal{O}(r^2(N^{d_x}+N^{d_v}))$ and $\mathcal{O}(r(N^{d_x} + N^{d_v}))$, respectively. The proposed method for a fixed rank has a somewhat larger overhead as the $3$ (in two dimensions) and $4$ (in three dimensions) moments have be integrated and stored as well. However, since for the proposed scheme in the $\varepsilon \to 0$ limit rank $1$ is sufficient, compared to a requirement for at least rank $6$/$10$ in two/three dimensions for the scheme in \cite{E18}, the proposed scheme is clearly more efficient. Note that computing a direct solution of the Boltzmann-BGK equation has at least a computational complexity of $\mathcal{O}(N^{d_x} N^{d_v})$, which is significantly more expensive.}
\end{remark}

\section{The fully discrete scheme \label{sec:fully-discrete}}

The scheme described in the previous section is still continuous in space. To provide an appropriate spatial discretization is the content of this section. Replacing spatial derivatives with their discrete counterpart has to be done in multiple places; namely, when computing $I_1^n$ and $I_2^n$ in equations \eqref{eq:I1n}-\eqref{eq:I2n}, when computing $K_j^{n+1}$ in equation \eqref{eq:Kn+1}, and when computing the coefficient $d_{ik}^1$ in equation \eqref{eq:d1}.

In this work we consider two approaches to spatial discretization. First, a method that uses FFT techniques (i.e., the Fourier spectral method) to approximate the derivatives is described in Section~\ref{subsec:fft}. It has the advantage that it is very accurate for smooth solutions and simple to implement. The second method, described in Section~\ref{subsec:cd}, is designed to handle solutions with sharp gradients/shock waves. This allows us to treat problems with small $\varepsilon$ where such features are common. This scheme is, at most, second order accurate.

\subsection{The Fourier spectral method \label{subsec:fft}}

For the Fourier spectral method we simply compute the spatial derivatives in equations \eqref{eq:I1n}, \eqref{eq:I2n}, \eqref{eq:Kn+1}, and \eqref{eq:d1} by an FFT. Let us emphasize that in place of the FFT a number of other spatial discretizations could be used. The choice of FFT here is mainly because for smooth solutions this results in a very accurate scheme and we further demonstrate that the dynamical low-rank approximation can be easily combined with a spectral discretization (which is not possible for the low-rank scheme described in \cite{Kormann15}, for instance).

\subsection{The shock capturing finite difference method \label{subsec:cd}}

    It is well-known that sharp gradients (for finite viscosity) and discontinuous solutions (for vanishing viscosity), so-called shock waves, are commonly encountered in fluid problems. Since our goal is to obtain a numerical scheme that can capture the fluid limit, we need an appropriate shock-capturing scheme. \textcolor{black}{Such a chock-capturing scheme has to be applied to both the $K_i$ in the low-rank approximation, see equation \eqref{eq:Kn+1}, and the moment equations given in \eqref{rho}-\eqref{rhou}. In the former the structure of the dynamical low-rank approximation allows us to relatively easily obtain the direction of flow and thus standard upwind or high resolution schemes can be employed. The latter is more delicate as the flux depends on the low-rank approximation and solving the corresponding Riemann problem would be rather challenging. We use the central scheme by Nessyahu and Tadmor to avoid having to obtain the direction of flow from the low-rank approximation. For simplicity, in this section we only consider the two-dimensional case} (i.e.~$d_x=2$). The extension to $d_x=3$ is immediate, \textcolor{black}{with the only major difference that the $K_i$ and the moments are three-dimensional.}

To discretize \eqref{eq:Kn+1} which is a symmetric hyperbolic system, we adopt the following approach. In the two dimensional setting, we have $c_{jl}^1 = [c_{jl}^{1;1}\, c_{jl}^{1;2}]^{\text{T}}$ and the matrices $c_{jl}^{1;m}$, for $m\in\{1,2\}$, are symmetric. Thus, there exist orthogonal matrices $T^{m}$ such that $\sum_{jl} T_{ij}^{m} c_{jl}^{1;m} T_{kl}^{m} = \lambda_{i}^{m} \delta_{ik}$. We apply the transformation $\hat{K}_i^n = \sum_{j} T_{ij}^1 K_j^n$ and obtain from equation \eqref{eq:Kn+1}
\[ \hat{K}_{i}^{n+1} =\frac{1}{1+\Delta t \rho^n/\varepsilon}\hat{K}_{i}^{n}-\frac{\Delta t}{1+\Delta t \rho^n/\varepsilon}\left[\lambda_i^1 \partial_{x}\hat{K}_{i}^{n} + \sum_{jl} T_{ij}^1 c_{jl}^{1;2} \partial_{y} K_l^n+\sum_{lj} T_{ij}^1 c_{jl}^{2}K_{l}^{n}\right]+\frac{\Delta t \rho^n}{\varepsilon+\Delta t \rho^n}\sum_j T_{ij}^1 \overline{V}_{j}. \]
Now, the direction of the flow regarding the advection in the $x$-direction is obvious and we can replace $\lambda_i^1 \partial_x \hat{K}_i^n$ by an appropriate discrete approximation, which we denote by $\delta_x(K_{\cdot}^n,\lambda_i^1)_i$. The simplest choice is to use upwinding for $\delta_x$. To obtain second order for smooth parts of the solution, we can use the Lax--Wendroff flux with the van Leer limiter as described in \cite[Chap. 16]{leveque1992}. Using an appropriate limiter is obviously important for problems that have sharp gradients/shocks. Finally transforming the system back to the original variables yields
\begin{equation} \label{eq:Kspace} K_{j}^{n+1} =\frac{1}{1+\Delta t \rho^n/\varepsilon}K_{j}^{n}-\frac{\Delta t}{1+\Delta t \rho^n/\varepsilon}\left[\sum_{i} T_{ij}^1 \delta_x(\hat{K}_{\cdot}^n,\lambda_i^1)_i + \sum_{l} c_{jl}^{1;2} \partial_{y} K_l^n+\sum_{l} c_{jl}^{2}K_{l}^{n}\right]+\frac{\Delta t \rho^n}{\varepsilon+\Delta t \rho^n} \overline{V}_{j}. \end{equation}
The discretization in the $y$-direction is treated similarly.

We now consider the update of $\rho$ and $u$ according to the system of conservation laws in equations \eqref{rho}-\eqref{rhou}. However, the flux depends on the Maxwellian $M$ and is thus a nonlinear (and rather complicated) function of $\rho$ and $u$. Obtaining the corresponding Riemann solver is, at best, a rather complicated undertaking. Another possible approach is to write down a scheme (which is typically nonlinear and involves a limiter for the shock-capturing purpose) at the level of the kinetic equation, i.e.~for equation \eqref{BGK}, and only after performing this spatial discretization the dynamical low-rank projection is applied. Such an approach has been used before, for example, to design AP kinetic schemes \cite{HJL17}. However, this does not work in the present setting as applying a nonlinear scheme at the kinetic level will result in a function which we can not compute using a convolution, as has to be done in {\bf Step 1} of the algorithm  (see Section~\ref{sec:first-order}) in order to obtain an efficient numerical method. 

To avoid the above issues, we choose to discretize \eqref{eq:I1n} and \eqref{eq:I2n}, which are the fluxes of the nonlinear conservation laws \eqref{rho}-\eqref{rhou}, using the central difference scheme by Nessyahu and Tadmor \cite{nessyahu1990}. This scheme avoids solving a Riemann problem and can be combined with a limiter to resolve sharp gradients without oscillations. In the following we will describe this method and how it can be applied to the present situation.
For a one-dimensional conservation law 
\[ \partial_t U(t,x) + \partial_x F(U(t,x)) = 0, \]
where $U$ is the vector of moments, the  Nessyahu--Tadmor scheme can be written in predictor-corrector form
\begin{align*}
    &U_j^{\star} = U_j^n - \tfrac{\Delta t}{2 \Delta x} F^{\prime}_j, \\
    &U_{j+1/2}^{n+1} = \tfrac{1}{2}(U_j^n +U_{j+1}^n) + \tfrac{1}{8}(U_j^{\prime}-U_{j+1}^{\prime}) - \tfrac{\Delta t}{\Delta x}(F(U_{j+1}^{\star})-F(U_{j}^{\star})).
\end{align*}
The choice of $F^\prime_j$ and $U^\prime_j$ is free as long as they satisfy $F^{\prime}_j/\Delta x = \partial_x F(U(t^n,x_j)) + \mathcal{O}(\Delta x)$ and $U_j^{\prime}/\Delta x = \partial_x U(t^n,x_j) + \mathcal{O}(\Delta x)$. We have to approximate these quantities in such a way that we obtain non-oscillatory solutions that preserve sharp gradients. As suggested in \cite{nessyahu1990} we use
\[ U_j^{\prime} = \text{MM}(U_{j+1}^n-U_{j}^n, U_{j}^n-U_{j-1}^n), \qquad 
   F_j^{\prime} = \text{MM}(F(U_{j+1}^n)-F(U_{j}^n), F(U_{j}^n)-F(U_{j-1}^n)),
\]
where $\text{MM}(x,y) $ is the commonly used minmod limiter. We note that this is a staggered scheme. That is, the output values $U^{n+1}_{j+1/2}$, are offset by half a grid point compared to the input values $U_j^n$.
In the present setting we use a two-dimensional generalization of this method given by (see, e.g., \cite{jiang1998})
\begin{align*} 
    &U_{ij}^{\star} = U_{ij}^n - \frac{\Delta t}{2 \Delta x} F_{ij}^{\prime x} - \frac{\Delta t}{2 \Delta y} G_{ij}^{\prime y}, \\
    &U_{i+1/2,j+1/2}^{n+1} = \frac{1}{4} (U_{ij}^n + U_{i+1,j}^n + U_{i,j+1}^n + U_{i+1,j+1}^n) \\
      & \qquad \qquad\qquad + \frac{1}{16} (U_{ij}^{\prime x} - U_{i+1,j}^{\prime x} + U_{i,j+1}^{\prime x} - U_{i+1,j+1}^{\prime x}) \\
      & \qquad \qquad\qquad + \frac{1}{16} (U_{ij}^{\prime y} - U_{i,j+1}^{\prime y} + U_{i+1,j}^{\prime y} - U_{i+1,j+1}^{\prime y}) \\
      & \qquad\qquad \qquad+ \frac{\Delta t}{\Delta x} (F(U_{i+1,j}^\star) - F(U_{ij}^\star)) + \frac{\Delta t}{\Delta y} (F(U_{i,j+1}^\star) - F(U_{ij}^\star)),
\end{align*}
where
\[ U_{ij}^{\prime x} = \text{MM}(U_{i+1,j}^n-U_{i,j}^n, U_{ij}^n-U_{i-1,j}^n), \qquad
   U_{ij}^{\prime y} = \text{MM}(U_{i,j+1}^n-U_{ij}^n, U_{ij}^n-U_{i,j-1}^n),
\]
and
\[ F_{ij}^{\prime x} = \text{MM}(F(U_{i+1,j}^n)-F(U_{i,j}^n), F(U_{ij}^n)-F(U_{i-1,j}^n)), \qquad
   F_{ij}^{\prime y} = \text{MM}(F(U_{i,j+1}^n)-F(U_{ij}^n), F(U_{ij}^n)-F(U_{i,j-1}^n)).
\]
We apply this scheme to \eqref{rho}-\eqref{rhou} which implies
\[ U=\left[\begin{array}{c}
\rho\\
\rho u
\end{array}\right],\qquad F(U)=\left[\begin{array}{c}
\sum_{ij}X_{i}^{n}S_{ij}^{n}\langle vV_{j}^{n}M^{n}\rangle_{v}\\
\sum_{ij}X_{i}^{n}S_{ij}^{n}\langle v\otimes vV_{j}^{n}M^{n}\rangle_{v}
\end{array}\right]. \]
The flux $F$ is computed using convolutions as described above.

While the Nessyahu--Tadmor scheme is also available in a non-staggered variant, we have used the staggered version here as we observed it to be more accurate in numerical simulations. Using a staggered scheme adds a small complication as we need the values at the original grid points to continue with the low-rank integrator. To accomplish this we do the following
\begin{enumerate}
    \item Apply the Nessyahu--Tadmor scheme with time step size $\Delta t/2$ and initial value $U_{ij}^n$ to obtain $U^{n+1/2}_{i+1/2,j+1/2}$.
    \item Obtain $X^n$ at the half grid points $(i+1/2,j+1/2)$ by computing averages between neighboring grid points. 
    \item Apply the Nessyahu--Tadmor scheme with time step size $\Delta t/2$ and initial value $U_{i+1/2,j+1/2}^{n+1/2}$ and the $X^n$ at the half grid points to obtain $U^{n+1}_{i,j}$. 
\end{enumerate}
This procedure takes the place of {\bf Step 1} in the algorithm described in Section~\ref{sec:first-order}.

The only ingredient that remains is the computation of the derivative in $d_{ik}^1$, for which we simply apply the one-sided first order finite difference. Since only the integral over such quantities enters in the computation, this was found to be sufficiently accurate. As an alternative a higher order stencil, which then would need to be limited in an appropriate way, could also be used here.

\section{Asymptotic analysis of the fully discrete scheme \label{sec:ap}}

In this section, we carry out a formal asymptotic analysis for the fully discretized scheme introduced above and show that for small values of $\varepsilon$, the scheme reduces to a consistent discretization of the NS equations.

First, our numerical scheme for the moment variables (updating $\rho^{n+1}$ and $u^{n+1}$) is equivalent to discretizing the equation (\ref{MM}) as follows
\begin{equation} \label{sch}
\frac{U^{n+1}-U^n}{\Delta t}+ \nabla_x\cdot  \langle v\phi M^ng^n\rangle_{v}=\mathcal{O}(\Delta x^p),
\end{equation}
where $p$ is the order of the spatial discretization for $\nabla_x\cdot  \langle v\phi M^ng^n\rangle_{v}$. 

Second, our numerical scheme for the low-rank variable (updating $g^{n+1}$ or $X_i^{n+1}$, $V_j^{n+1}$ and $S_{ij}^{n+1}$) is equivalent to discretizing the equation (\ref{hh}) as
\begin{equation} \label{gg-scheme}
\frac{g^{n+1}-g^n}{\Delta t}=-v\cdot\nabla_{x}g^n-\mathcal{M}(U^n)g^n+\frac{\nu}{\varepsilon}(1-g^{n+1})+\mathcal{O}(\Delta t+\Delta x^p),
\end{equation}
provided the same order of spatial discretization is used for $v\cdot\nabla_{x}g^n$. Note that $\mathcal{M}(U)$ in the case of $T=1$ is given in (\ref{MM123}). The $\mathcal{O}(\Delta t)$ term is due to the fact that the projector splitting integrator has an additional first order time error.

Rearranging (\ref{gg-scheme}), we have
\begin{align} \label{schg}
g^{n+1}&=1-\frac{\varepsilon}{\nu}\left[\frac{g^{n+1}-g^n}{\Delta t} +v\cdot\nabla_{x}g^n+\mathcal{M}(U^n)g^n\right]+\mathcal{O}(\varepsilon\Delta t+\varepsilon\Delta x^p)\nonumber\\
&=1-\frac{\varepsilon}{\nu}\mathcal{M}(U^n)+\mathcal{O}(\varepsilon^2+\varepsilon\Delta t+\varepsilon\Delta x^p)\nonumber\\
&=1-\frac{\varepsilon}{\nu}\mathcal{M}(U^{n+1})+\mathcal{O}(\varepsilon^2+\varepsilon\Delta t+\varepsilon\Delta x^p),
\end{align}
where we used $g=1+\mathcal{O}(\varepsilon)$ (the first equality implies so) to derive the second line. Finally, substituting (\ref{schg}) into (\ref{sch}) yields
\begin{equation} 
\frac{U^{n+1}-U^n}{\Delta t}+\nabla_x\cdot \langle v\phi M^{n}\rangle_{v}=\frac{\varepsilon}{\nu}\nabla_x\cdot \langle v\phi M^n\mathcal{M}(U^{n})\rangle_{v}+\mathcal{O}(\Delta x^p+\varepsilon\Delta t+\varepsilon\Delta x^p+\varepsilon^2),
\end{equation}
which is a consistent first order in time and $p$-th order in space discretization to the NS equations (\ref{NS}).

\section{Numerical results \label{sec:numerical}}

In this section we will consider three test problems. Namely, a shear flow, an explosion, and the
relaxation of a plasma beam. The first two examples operate close to the fluid regime and thus show
that our numerical method can capture that limit. The explosion, in particular, demonstrates that
using the Nessyahu--Tadmor central scheme we can handle shock waves. The beam relaxation problem, on
the other hand, is run for larger $\varepsilon$ and is thus fully kinetic. In this context we will also
demonstrate that our dynamical low-rank scheme works well for spatially varying $\varepsilon$; i.e.~in situations where $\varepsilon$ can become extremely small in some regions and large in others. 

All simulations in this section are conducted in $2+2$ dimensions (i.e.~$d_x=2$ and $d_v=2$) and periodic boundary conditions are imposed in all spatial and velocity directions. For the spatial discretization, we will, depending on the flow regime, use either the Fourier spectral method (see Section~\ref{subsec:fft}) or the finite difference method (see Section~\ref{subsec:cd}). The latter will be referred to as SCFD (shock-capturing finite differences) in the following. 

\subsection{Shear flow}

Shear flow is a classical example in computational fluid dynamics. A common setup (see, e.g., \cite{bell1989second,liu2000high,einkemmer2014}) is
\begin{align} \label{eq:shear-flow-iv}
    \begin{split}
    \rho(0,x,y) &=1,\\
    u_1(0,x,y) &=v_{0}\begin{cases}
    \tanh\left(\frac{y-\tfrac{1}{4}}{\Delta}\right), & y\leq\tfrac{1}{2},\\
    \tanh\left(\frac{\tfrac{3}{4}-y}{\Delta}\right), & y>\tfrac{1}{2},
    \end{cases}\\
    u_2(0,x,y) &=\delta\sin(2\pi x),
    \end{split}
\end{align}
where $(x,y) \in [0,1]^2$, $v_0=0.1$, $\Delta = 1/30$, and $\delta=5 \cdot 10^{-3}$. Physically the problem models a fluid that moves to the left at the top and bottom of the domain and to the right in the middle of the domain. The small perturbation in the velocity then grows in time and creates a dynamics dominated by vortex structures. 

For our dynamical low-rank integrator the moments $\rho$, $u_1$, and $u_2$, as specified above, are used to initialize a Maxwellian according to equation \eqref{Max} in the domain $(v_1,v_2) \in [-6,6]^2$. The Reynolds number Re is prescribed and then $\varepsilon$ is chosen according to $\varepsilon=v_{0}/\text{Re}$. That is, we use the flow velocity $v_0$ as the characteristic velocity to define the Reynolds number. Note that this is not entirely consistent with Remark \ref{rem:reynolds}, where the speed of sound is used instead. However, since the former is the usual choice in the weakly compressible regime, this makes comparison with the literature (e.g.~\cite{bell1989second,liu2000high,einkemmer2014}) easier.

\begin{figure}[h!]
    \centering Fluid solver

    \includegraphics[width=7cm]{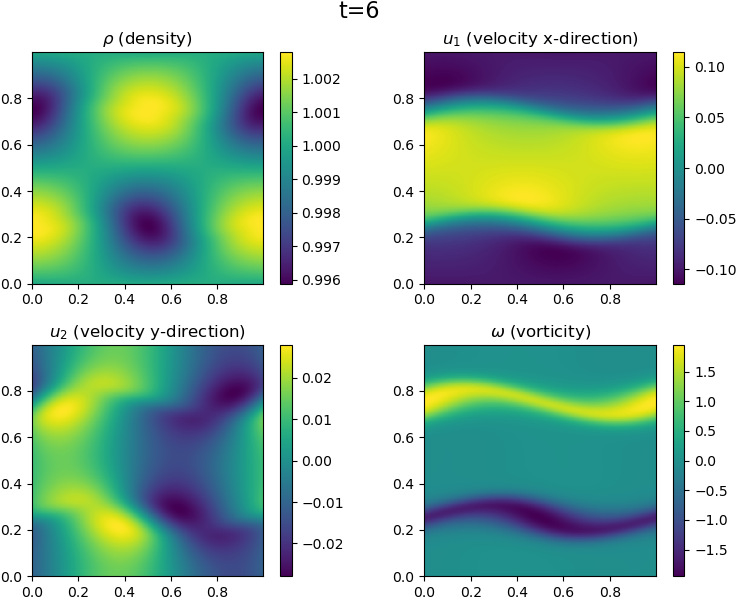}
    \includegraphics[width=7cm]{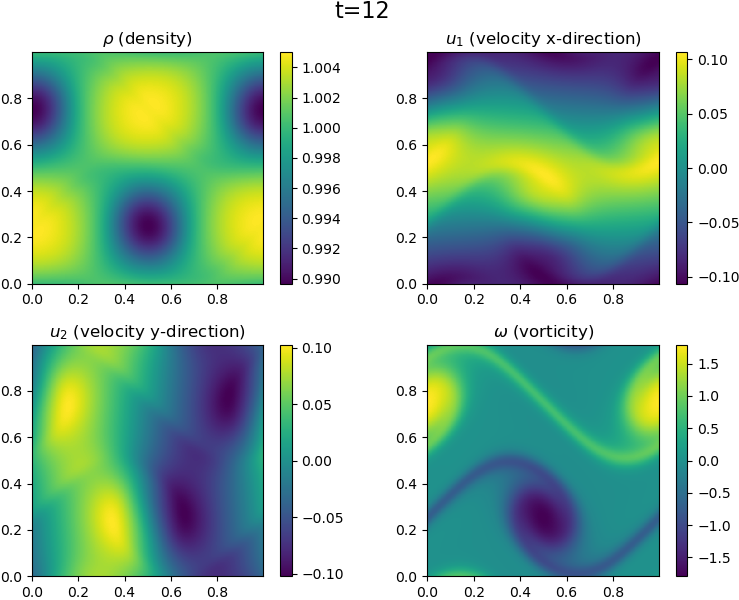}

    \vspace{0.5cm}
    Dynamical low-rank (Fourier spectral method)

      \includegraphics[width=7cm]{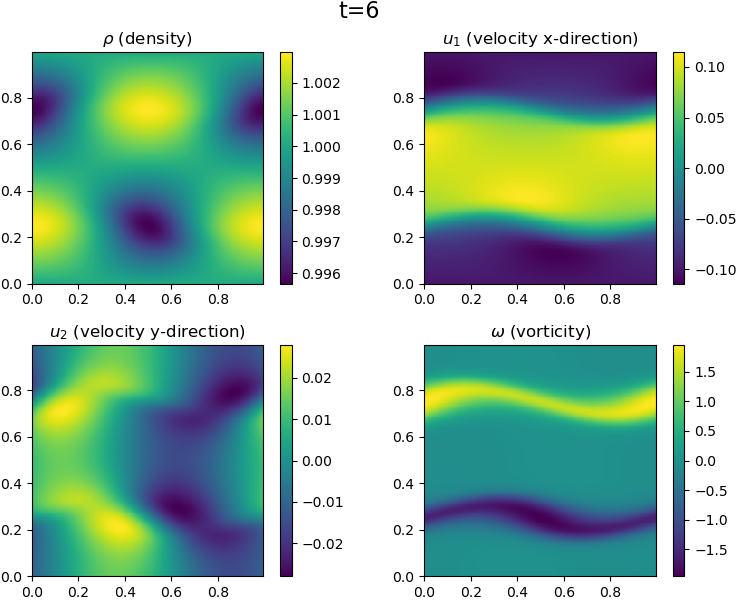}
    \includegraphics[width=7cm]{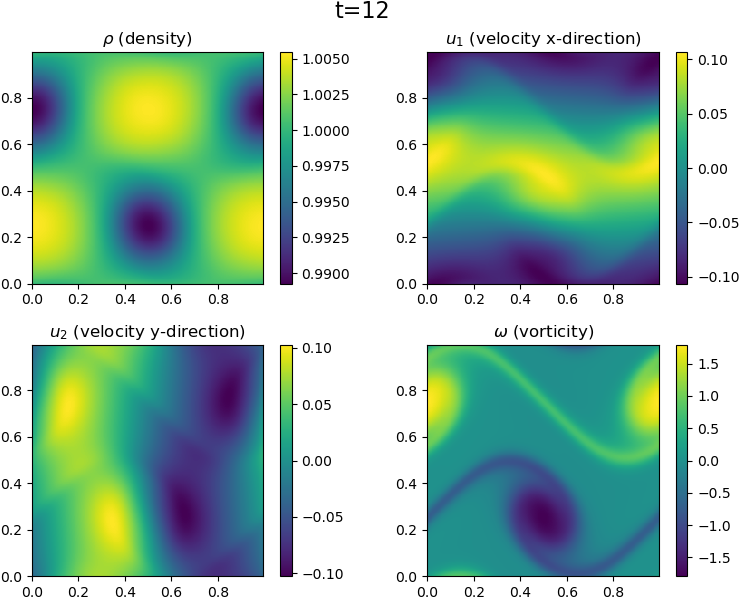}

    \caption{Numerical simulation of the shear flow \eqref{eq:shear-flow-iv} at time $t=6$ (left) and $t=12$ (right) with $\text{Re}=10^3$ is shown. The fluid solver (top) uses a MacCormack method, $512$ grid points in each direction, and is run using a CFL number of 0.9 ($\Delta t \approx 1.5 \cdot 10^{-3}$). The dynamical low-rank integrator (bottom) uses rank $r=3$ and $256$ grid points in each spatial direction and $32$ grid points in each velocity direction and time step size $\Delta t = 2 \cdot 10^{-4}$. \label{fig:shear-flow-Re1000}}
\end{figure}

\begin{figure}[h!]
    \centering
    \includegraphics[width=7cm]{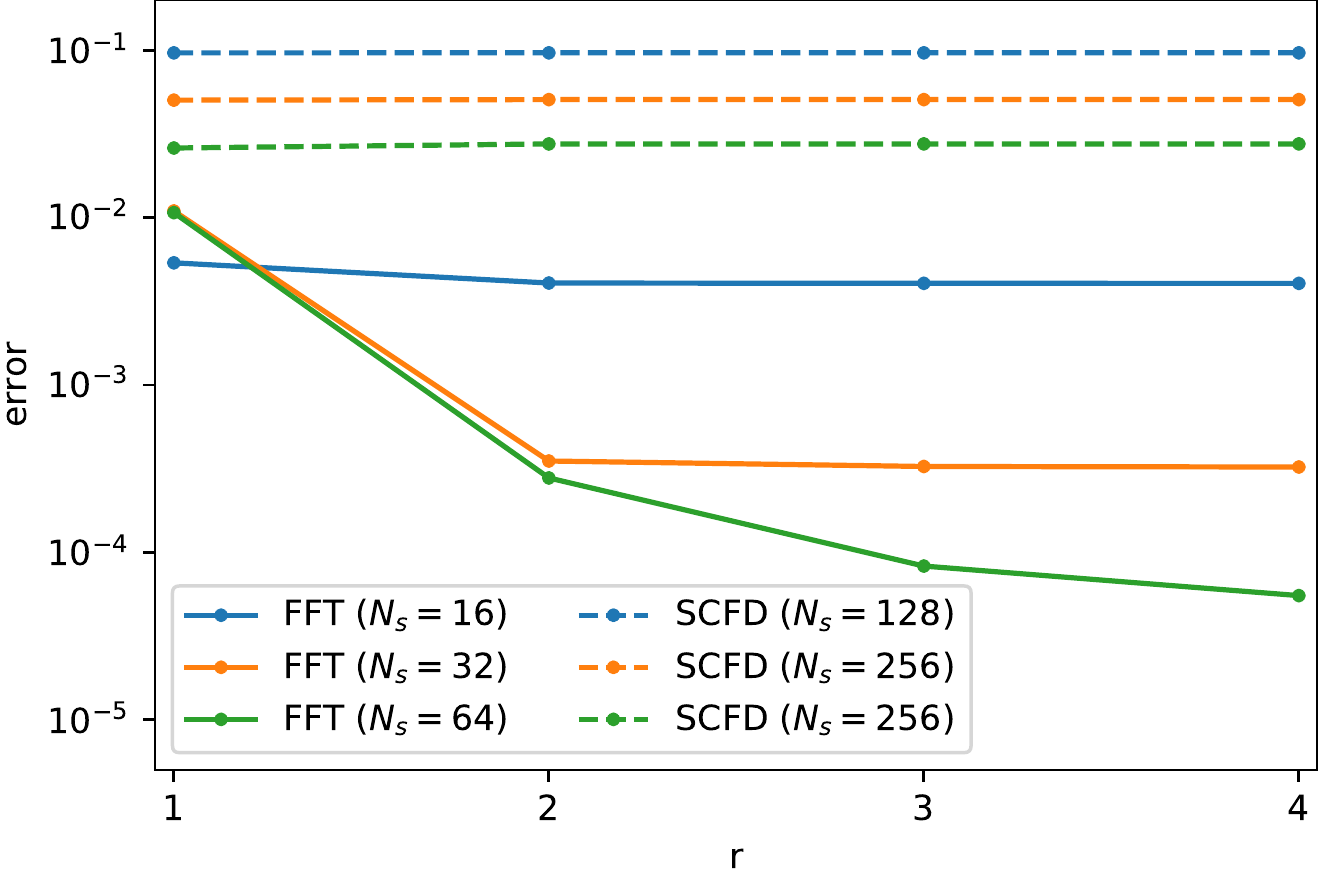}
    \caption{The error in the moments, i.e.~$\max\{\Vert \rho - \rho_{\text{ref}} \Vert_{\infty}, \Vert \rho u - \rho u_{\text{ref}} \Vert_{\infty} \}$, at time $t=2$ is shown as a function of the number of grid points in the $x$ and $y$ directions (denoted by $N_s$ in the plot) and the rank $r$. In all cases $32$ grid points in the velocity direction, a time step size of $\Delta t = 2\cdot 10^{-5}$, and $\text{Re}=1000$ is used. Results are shown both for the SCFD and the FFT based method.  We note that the error due to the low-rank approximation is below $10^{-2}$ for $r=1$ and below $10^{-4}$ for $r=3$. The reference solution $(\rho_{\text{ref}}, \rho u_{\text{ref}})$ is computed using a MacCormack based fluid solver with $2048 \times 2048$ grid points. \label{fig:quantitative}}
\end{figure}

In Figure \ref{fig:shear-flow-Re1000} a comparison at $\text{Re}=1000$ between a fluid solver based on the MacCormack method \cite{maccormack} and the proposed dynamical low-rank integrator is conducted. We compare density $\rho$, velocity $u_1$, $u_2$, and vorticity $\omega=\partial_x u_2 - \partial_y u_1$ and observe excellent agreement with the fluid solver, as expected. \textcolor{black}{A more quantitative comparison is conducted in Figure \ref{fig:quantitative}}.

\begin{figure}[h!]
    \centering
    \includegraphics[width=7cm]{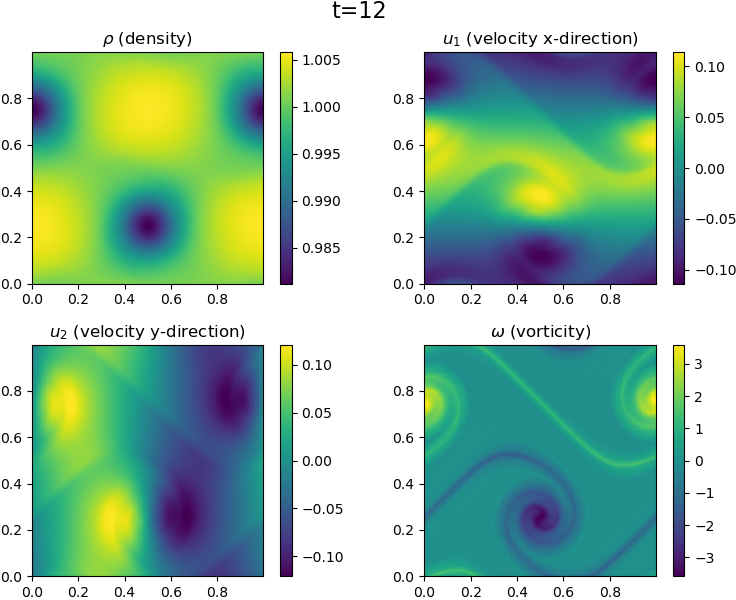}
    \includegraphics[width=7cm]{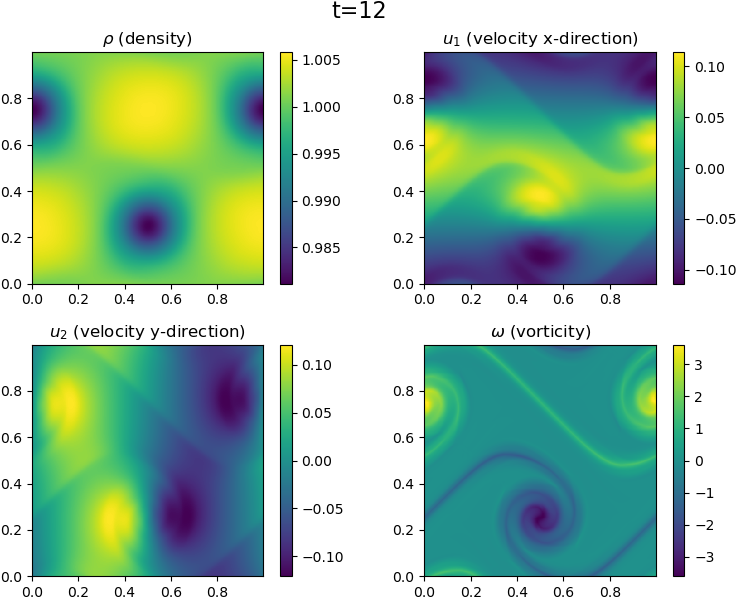}

    \caption{Numerical simulation of the shear flow \eqref{eq:shear-flow-iv} at time $t=12$, $\text{Re}=10^5$, and with rank $r=2$ (left) and $r=10$ (right) is shown. The SCFD implementation with $512$ grid points in each spatial and $32$ grid points in each velocity direction and a time step size of $\Delta t=10^{-3}$ is used.\label{fig:shear-flow-Re1e5}}
\end{figure}

Next, we demonstrate that increasing the Reynolds number (i.e.~decreasing $\varepsilon$) presents no difficulty for our approach. In this case we have to use the SCFD implementation as the Fourier spectral method, as one would expect, produces spurious oscillations and shows stability issues in this regime. The corresponding results are shown in Figure \ref{fig:shear-flow-Re1e5} and agree well with the results that have been reported in the literature. 

We also remark that choosing a small rank ($r=2$) is sufficient to obtain very accurate results. Let us emphasize that this is in contrast to earlier work \cite{E18}, where the expansion of the equilibrium mandates that we use a certain rank (at least $6$ in 2D and $10$ in 3D). In the method proposed here no assumption on the low-rank structure of the equilibrium needs to be made and thus even a very small rank is sufficient to accurately represent problems which are well described by the fluid limit. \textcolor{black}{We also note that in this limit, as has already been observed in \cite{E18}, a coarse resolution in the velocity directions is sufficient.} This is precisely what we observe in Figure \ref{fig:shear-flow-Re1e5}.

\begin{figure}[h!]
    \centering
    $Re=500$

    \includegraphics[width=7cm]{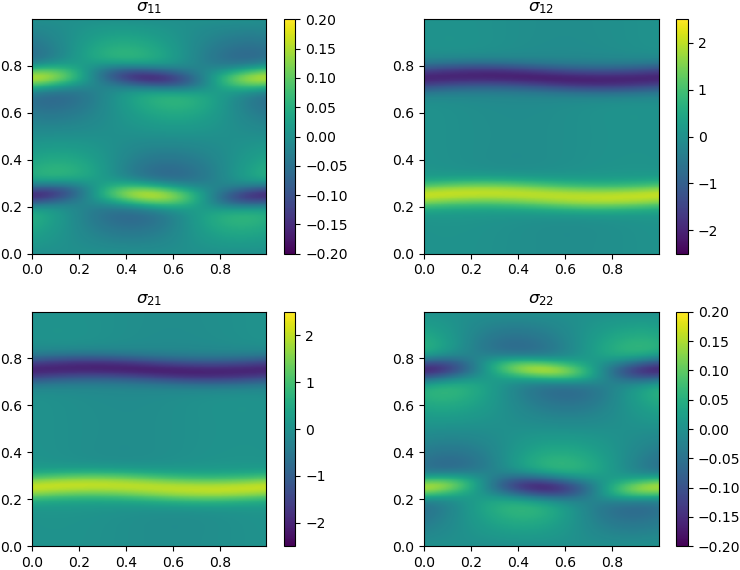}
    \includegraphics[width=7cm]{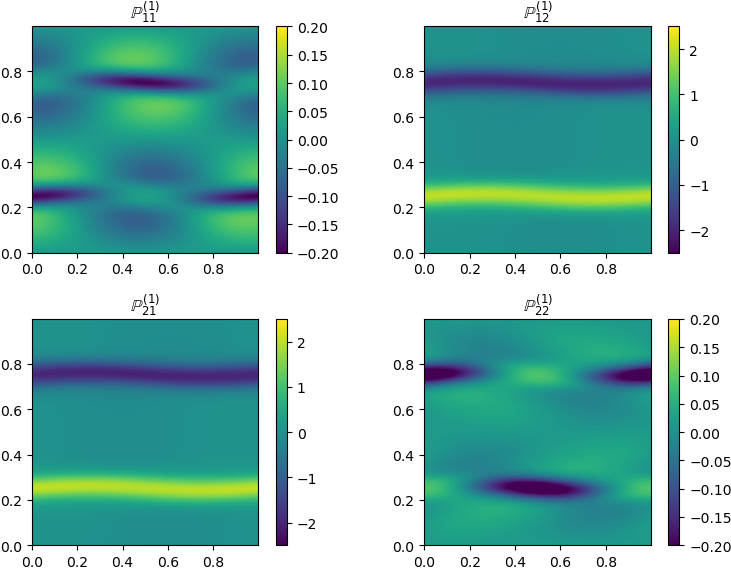}
    
    $Re=4000$
    
    \includegraphics[width=7cm]{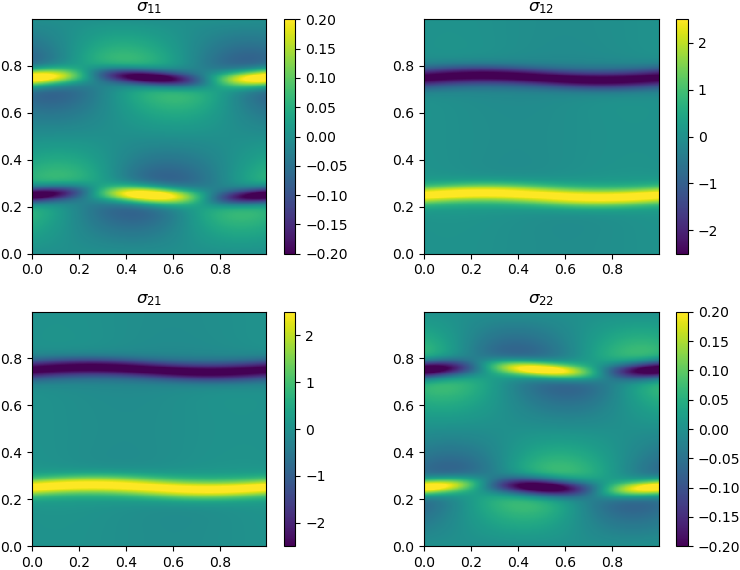}
    \includegraphics[width=7cm]{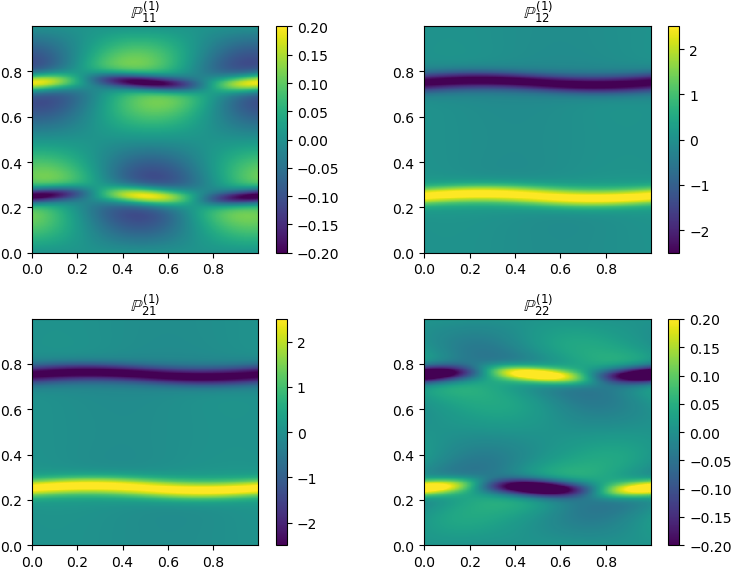}

    \caption{Comparison of $\sigma(u)$, computed with the MacCormack based fluid solver, with the stress tensor $\mathbb{P}^{(1)}$, computed with the dynamical low-rank approximation, for $Re=500$ and $Re=4000$. All simulations are conducted with rank $r=5$, $512$ grid points in each spatial direction, $64$ grid points in each velocity direction, and a time step size of $\Delta t=5\cdot 10^{-5}$. \label{eq:stresstensor1}}
\end{figure}

\begin{figure}[h!]
    \centering
    \includegraphics[width=7cm]{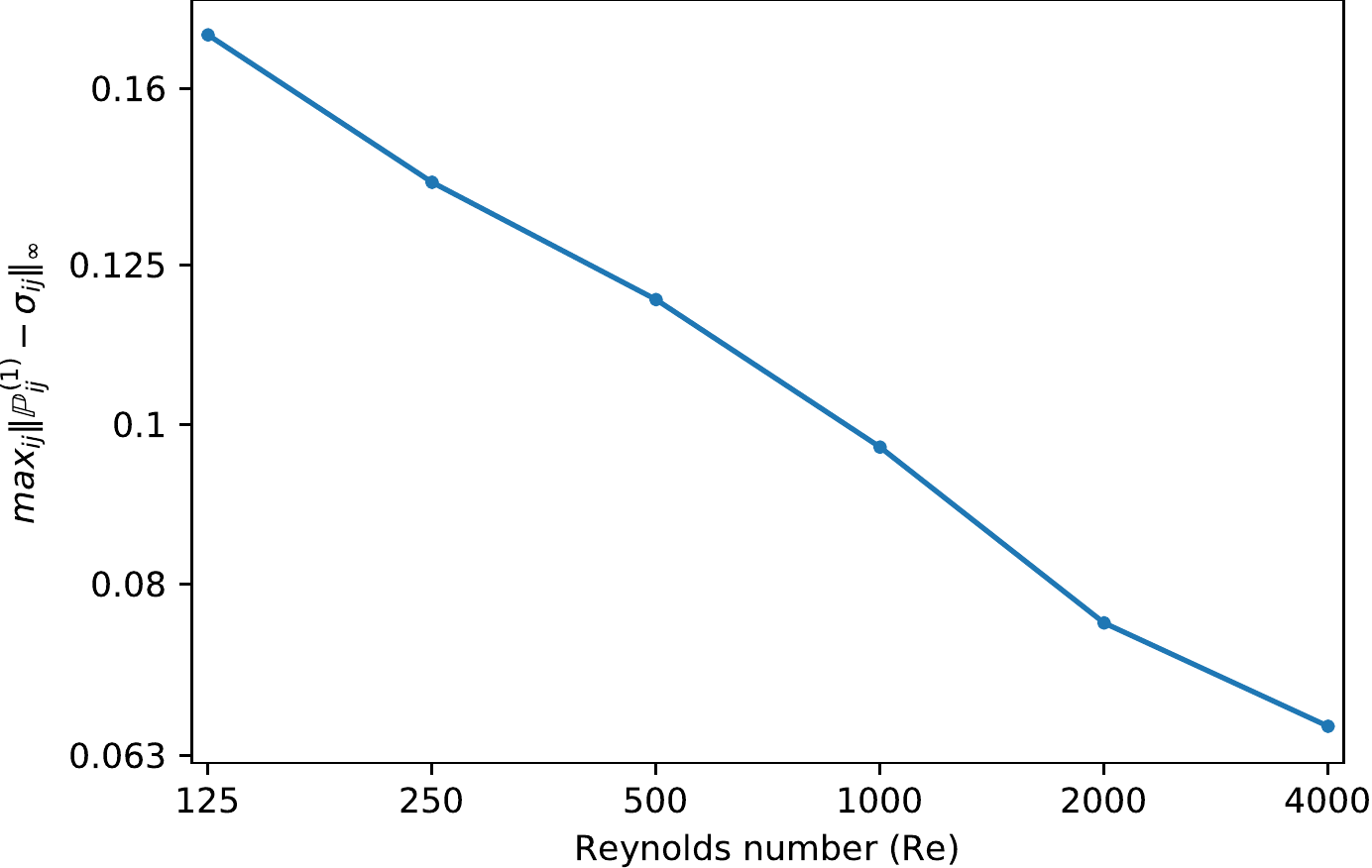}
    \caption{Difference between $\sigma(u)$, computed with the fluid solver, and the stress tensor $\mathbb{P}^{(1)}$, computed with the dynamical low-rank approximation, in the maximum norm and as a function of the Reynolds number is shown. All simulations are conducted until $t=2$ and with $r=5$, $512$ grid points in each spatial direction, $64$ grid points in each velocity direction, and a time step size of $\Delta t=5\cdot 10^{-5}$. \label{eq:stresstensor2} }
\end{figure}

Finally, we investigate the convergence of the stress tensor computed by the kinetic integrator (i.e.~the proposed dynamical low-rank algorithm) with the corresponding term in the NS equation. That is, we compare
\begin{equation} \mathbb{P}^{(1)} = -\frac{1}{\varepsilon} \int(v-u)\otimes(v-u)(f-M)\,dv 
= \frac{1}{\varepsilon}\bigg[-\frac{\rho}{(2\pi)^{d_v/2}}\sum_{ij}X_{i}S_{ij}\int(v-u)\otimes(v-u)e^{-\vert v-u\vert^{2}/2}V_{j}(v)\,dv+\rho \text{Id}\bigg], \label{eq:stresstensor} \end{equation}
obtained by the dynamical low-rank integrator, with $\sigma(u)$ (defined in (\ref{sigmau})) obtained from the fluid solver (here $d_v=2$). To evaluate equation \eqref{eq:stresstensor}, without incurring the cost of forming the entire density function, we once again write the integral as a convolution (as described in Section~\ref{sec:first-order}). We know from theory that $\mathbb{P}^{(1)}$ converges to $\sigma(u)$ as $\varepsilon\to 0$ (see Section~\ref{sec:fluid-limit}). This is also observed in the numerical results shown in Figure \ref{eq:stresstensor1} and Figure \ref{eq:stresstensor2}. In Figure \ref{eq:stresstensor2} we have plotted the difference between these two quantities as a function of the Reynolds number.

\subsection{Explosion}

In this section we consider an initial overpressure in a small region that expands and results in a propagating shock wave, i.e.~an explosion. We will use the following initial value
\begin{align} \label{eq:explosion}
    \begin{split}
  \rho(0,x,y) &=\begin{cases}
    1, & x^{2}+y^{2}\leq R^{2},\\
    \tfrac{1}{10}, & \text{otherwise},
    \end{cases} \\
        u_1(0,x,y) &=0, \qquad u_2(0,x,y) =0,
    \end{split}
\end{align}
with $(x,y) \in [-1.5,1.5]^2$ and $R=10^{-2}$. As before, the specified moments are used to initialize a Maxwellian according to equation \eqref{Max} in the domain $(v_1,v_2) \in [-6,6]^2$. This is a good example to demonstrate that our scheme is able to resolve sharp gradients/shocks in a compressible flow. 

\begin{figure}[h!]
    \centering
    \includegraphics[width=7cm]{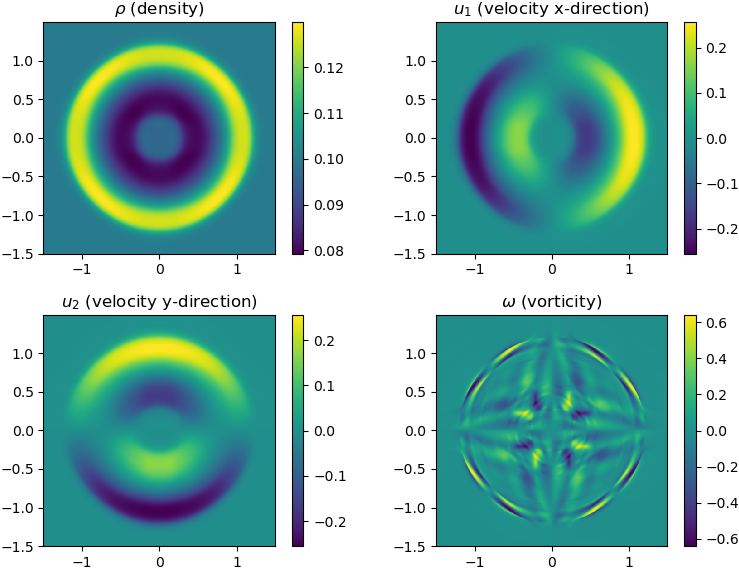}
    \hfill
    \includegraphics[width=7cm]{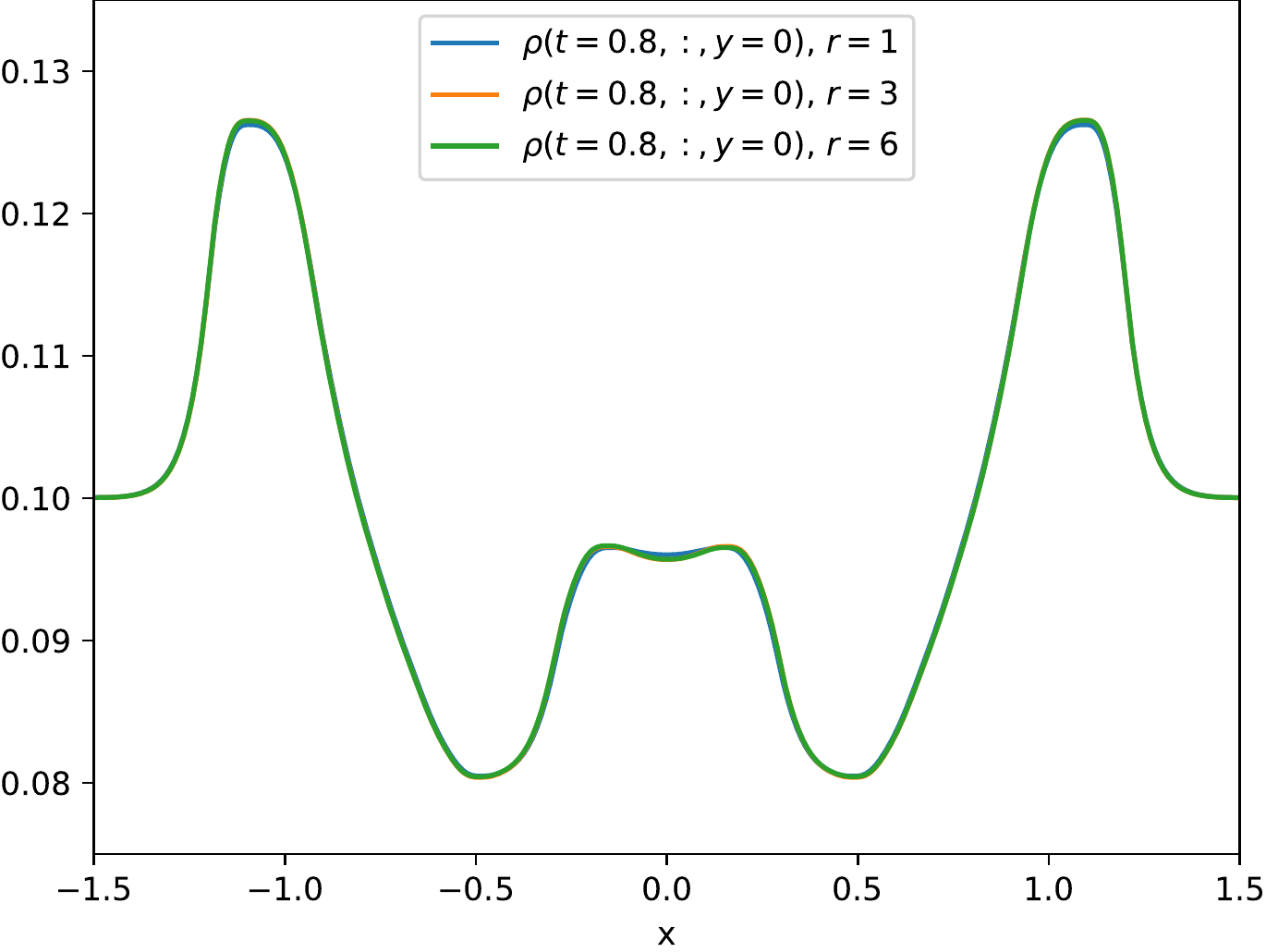}

    \caption{Numerical simulation of the explosion \eqref{eq:explosion} at time $t=0.8$, $\varepsilon=10^{-5}$, and with rank $r=6$ is shown on the left. The SCFD implementation with $512$ grid points in each spatial and $32$ grid points in each velocity direction and a time step size of $\Delta t=5 \cdot 10^{-5}$ is used. On the right a slice through the solution, at $y=0$, is shown for $r=1$, $r=3$, and $r=6$. \label{fig:explosion}}
\end{figure}

The numerical results are shown in Figure \ref{fig:explosion}. We clearly observe the shock front that radially propagates outward starting from the center. Let us also emphasize that, for $\varepsilon = 10^{-5}$, a small rank is sufficient to resolve the dynamics; in fact, starting from $r=3$ the solutions produced by the low-rank integrator are indistinguishable in the plot. %

\subsection{Beam relaxation}

In this section we consider a fully kinetic problem, i.e.~a problem where the dynamics is not well represented by the fluid limit. A beam with a constant velocity and Maxwellian distribution propagates into a stationary gas. The initial value is given by \cite{crouseilles2005hybrid}
\begin{equation} \label{eq:beam}
    f(0,x,y,v,w) = \frac{1}{\sqrt{2\pi}}\left(e^{-(v^{2}+w^{2})/2}+n_{b}\mathrm{e}^{-((v-v_{b})^{2}+(w-w_{b})^{2})/(2T_{b})}\right),
\end{equation}
where $(x,y,v,w) \in [0,1]^2 \times [-8,8]^2$, $n_{b}=10^{-3}$, $v_{b}=4$ , $w_{b}=2$, $T_{b}=0.1$. This is a rank 1 initial value that in our framework, i.e. with $f=Mg$, is represented as follows
\[ \rho = 1, \qquad u_1 = u_2 = 0,\]
and
\[ g=1+n_{b}\mathrm{e}^{-((v-v_{b})^{2}+(w-w_{b})^{2})/(2T_{b})+(v^{2}+w^{2})/2}. \]

The expected behavior is that the perturbation (the beam) in velocity space diminishes over time according to $\exp(-t/\varepsilon)$. That is, the collisions thermalize the beam until it is indistinguishable from the background gas. This problem has been considered in the context of hybrid kinetic-fluid models in \cite{crouseilles2005hybrid}. The hybrid method proposed there, which treats fast particles kinetically and the slow particles by a fluid model, has difficulties in capturing the physical behavior. 

\begin{figure}
    \centering
    \includegraphics[width=7cm]{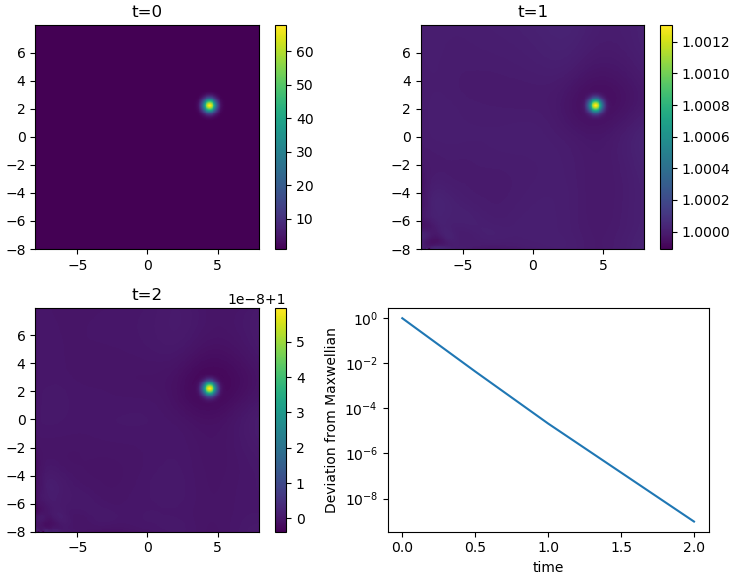}
    \hfill
    \includegraphics[width=7cm]{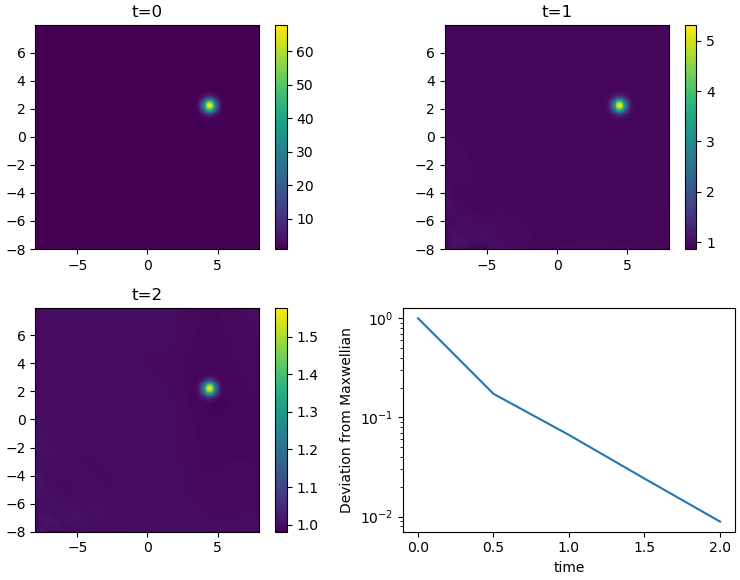}

    \caption{We plot $g(0,0,:,:)$ for the beam relaxation problem \eqref{eq:beam} at time $t=0$, $t=1$, and $t=2$ and for $\varepsilon=0.1$ (left) and $\varepsilon=0.5$ (right). The SCFD implementation with $r=10$, $32$ grid points in each spatial and $256$ grid points in each velocity direction, and a time step size of $\Delta t=10^{-4}$ is used. On the bottom right the maximal deviation from equilibrium, i.e.~the deviation in the maximum norm (normalized to the initial value), is shown as a function of time. \label{fig:beam}}
\end{figure}

The numerical results using the low-rank integrator are presented in Figure \ref{fig:beam}. \textcolor{black}{Note that because the initial value is homogeneous in the $x$ directions we can use a relatively coarse resolution in space.} For $\varepsilon=0.1$/$\varepsilon=0.5$ the deviation from equilibrium reduces by approximately $10^{-9}$/$10^{-2}$, which clearly demonstrates that the low-rank algorithm accurately reproduce the expected decay rates.  %

\subsection{Beam relaxation with spatially varying $\varepsilon$}

One interesting application of the proposed dynamical low-rank integrator is in applications where the Knudsen number $\varepsilon$ varies over several orders of magnitude (see the discussion on hybrid kinetic/fluid schemes in the Introduction). We explore this here in the context of the beam relaxation problem introduced in the previous section with a spatially varying $\varepsilon$ given by
\begin{equation} \label{eq:vareps}
    \varepsilon(x,y)=\varepsilon_{0}+\tanh(1-11x)+\tanh(1+11x), 
\end{equation}
with $\varepsilon_{0}=10^{-4}$ in the spatial domain $(x,y)\in[-1,1]\times[0,1]$. Thus, we have $\varepsilon \approx 1$ in the middle of the domain (kinetic regime) and $\varepsilon \approx 10^{-4}$ (fluid regime) as we get closer to the boundary. 

\begin{figure}
    \centering
    \begin{center}\textbf{Dynamical low-rank with $r=3$}\end{center}
    \includegraphics[width=12cm]{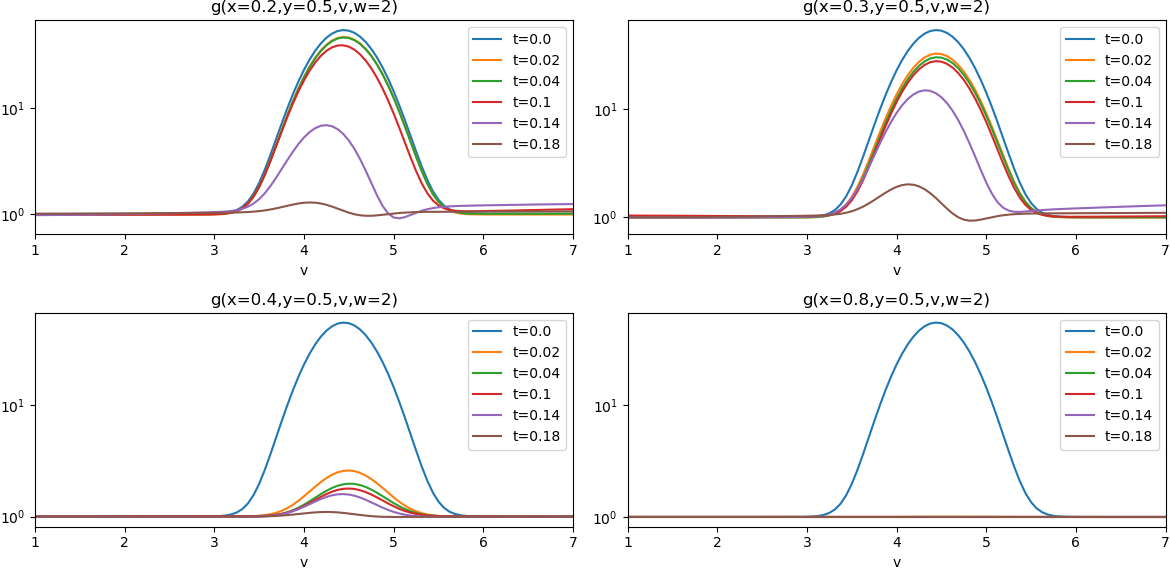}

    \begin{center}\textbf{Dynamical low-rank with $r=10$}\end{center}
    \includegraphics[width=12cm]{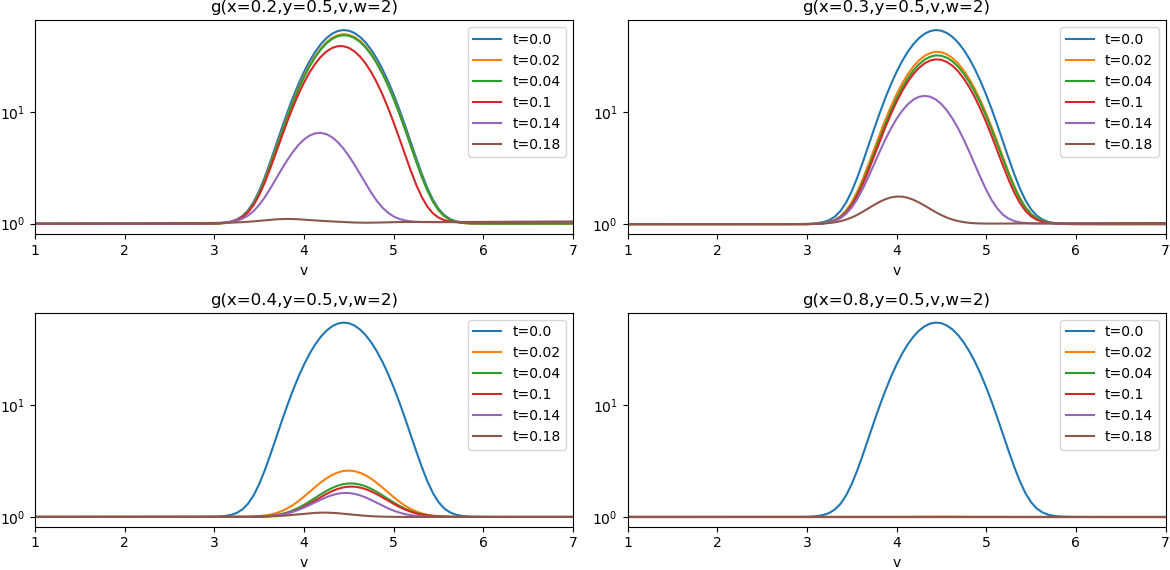}

    \begin{center}\textbf{Dynamical low-rank with $r=60$}\end{center}
    \includegraphics[width=12cm]{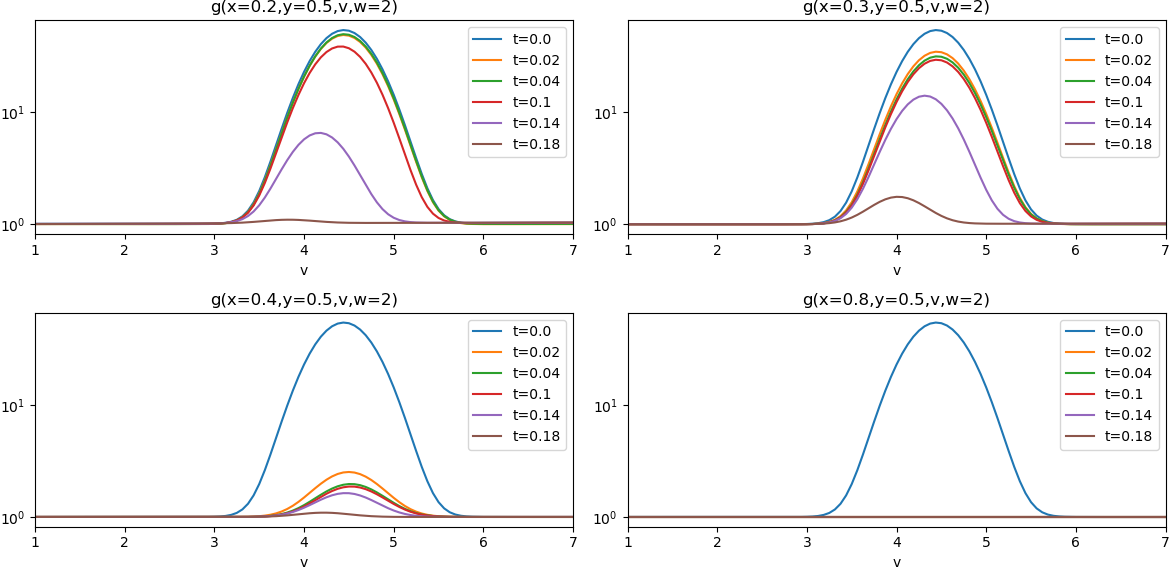}
    \caption{We plot slices through $g$ corresponding to $x=0.2, 0.3, 0,4, 0.8$ and $y=0.5$, $w=2$ for the beam relaxation problem \eqref{eq:beam} with spatially varying $\varepsilon$ given by equation \eqref{eq:vareps}. In all directions (both space and velocity) $256$ grid points are employed and the SCFD scheme (with $\delta_x$ in equation \eqref{eq:Kspace} discretized by upwinding) and a time step size of $\Delta t = 5 \cdot 10^{-5}$ is used in all simulations. \label{fig:vareps} }
\end{figure}

The numerical results obtained are shown in Figure \ref{fig:vareps}. We investigate here the deviation from equilibrium as a function of velocity. We observe that the thermalization proceeds with different speeds depending on the value of $\varepsilon$ on that part of the domain, as we would expected. We have also compared numerical simulations obtained with different ranks. For $r=3$ a significant deviation from the reference solution (computed with $r=60$) can be observed. However, starting with rank $r=10$ the low-rank error is small and we observe excellent agreement with the reference solution.

\bibliographystyle{plain}
\bibliography{hu_bibtex,le_bibtex}
\end{document}